\documentclass[11pt,reqno]{amsart}
\usepackage[T1]{fontenc}
\usepackage{lmodern,mathtools,amssymb,microtype}
\usepackage[letterpaper,hmargin=1.1in,top=0.85in,bottom=0.86in]{geometry}
\usepackage[hidelinks]{hyperref}
\usepackage{comment}
\hypersetup{pdftitle={Multiplication from comultiplication for canonical bases of modified quantum groups},pdfauthor={Jiepeng Fang and Xuhua He}}
\newtheorem{theorem}{Theorem}[section]
\newtheorem{lemma}[theorem]{Lemma}
\newtheorem{remark}[theorem]{Remark}
\newtheorem{proposition}[theorem]{Proposition}
\newtheorem{corollary}[theorem]{Corollary}
\newcommand{\U}{\mathbf U}
\newcommand{\B}{\mathbf B}

\numberwithin{equation}{section}

\usepackage{xcolor}
\usepackage[normalem]{ulem}

\title[Multiplication from comultiplication]{Multiplication from comultiplication in modified quantum groups}

\author{Jiepeng Fang}
\address[J. Fang]{Department of Mathematics and New Cornerstone Science Laboratory,
The University of Hong Kong, Pokfulam, Hong Kong SAR, China}
\email{fangjp@hku.hk}
\author{Xuhua He}
\address[X. He]{Department of Mathematics and New Cornerstone Science Laboratory,
The University of Hong Kong, Pokfulam, Hong Kong SAR, China}
\email{xuhuahe@hku.hk}

\date{}
\begin{document}

\begin{abstract}
We establish an identity relating the multiplication in a modified quantum group with the comultiplication in a larger one, both expressed in Lusztig's canonical bases. In symmetric type, the coefficients in the product of two canonical basis elements of the modified quantum group occur, term for term, in the coproduct of a single canonical basis element of the larger one. Combined with our earlier positivity theorem for comultiplication, this proves Lusztig's conjecture that the structure constants of the multiplication are Laurent polynomials with nonnegative integer coefficients for every symmetric Cartan datum.
\end{abstract}

\maketitle

\section*{Introduction}

Quantum groups were introduced independently by Drinfeld and Jimbo in the mid-1980s in their work on the quantum Yang--Baxter equation \cite{Drinfeld,Jimbo}. Lusztig introduced the modified form $\dot\U$ of a quantum group, in which a system of mutually orthogonal idempotents replace the Cartan subalgebra \cite[Chapter~23]{Lusztig-1993}. The theory of canonical bases was introduced by Lusztig \cite{Lusztig-1990,Lusztig-1991,Lusztig-1993}. 

For the positive and negative parts of a quantum group associated with a symmetric Cartan datum, Lusztig proved that the structure constants of multiplication and comultiplication with respect to their canonical bases are Laurent polynomials in $v$ with nonnegative integer coefficients. These positivities come from Lusztig's geometric construction. The canonical basis elements correspond to simple perverse sheaves, and the multiplication corresponds to the convolution. The convolution of two such sheaves decomposes as a direct sum of shifts of simple perverse sheaves, and the coefficients record their multiplicities in each degree \cite[Theorem 14.4.13]{Lusztig-1993}.

The canonical basis $\dot{\mathbf{B}}$ of the modified quantum group $\dot\U$ involves both the canonical bases of the positive and negative parts. Positivity
for the multiplication in each part separately therefore does not by itself give the positivity for the multiplication in $\dot\U$. In 1993,
Lusztig conjectured that, for every symmetric Cartan datum, the structure constants of the multiplication and comultiplication in $\dot\U$ with respect to the canonical basis $\dot\B$ belong to $\mathbb N[v,v^{-1}]$ \cite[Conjecture~25.4.2]{Lusztig-1993}.

In our earlier work, we proved the comultiplication assertion for
every symmetric Cartan datum \cite[Theorem~6.3]{Fang-He-2026}.
In this paper, we deduce the multiplication assertion from that
result. The key is an identity: the coefficients in the product 
of two given canonical basis elements occur, term for term, in the coproduct of a single canonical basis element of a larger modified quantum group $\dot{\tilde{\mathbf{U}}}$. The larger $\dot{\tilde{\mathbf{U}}}$ is obtained by the \emph{thickening} construction of \cite{Fang-He-2026} which attaches one new vertex by a single edge to each vertex of the original Dynkin-Coxeter diagram.

\begin{theorem}[{\textbf{=Theorem \ref{thm:main}, Proposition \ref{prop:product}}}]\label{thm:intro-comparison}
Assume that the root datum is simply connected and the Cartan datum
is symmetric. Let $\dot{b}\in \dot{\mathbf{B}}$ and $\dot{b}'\in \dot{\mathbf{B}}\cap \dot{\mathbf{U}}1_{\zeta}$ be canonical basis elements of the modified quantum group $\dot{\mathbf{U}}$ such that $\dot{b}\dot{b}'\neq0$. Then for a suitable lift $\tilde{\zeta}'$ of $\zeta'$ and any sufficiently large $N\in \mathbb{N}$, there are canonical basis elements $B_N,H_N\in \dot{\tilde{\mathbf{B}}}$ of the larger modified quantum group $\dot{\tilde{\mathbf{U}}}$ such that
\[
 \dot{\Delta}_{\tilde\zeta',0}(B_N)
   = \iota(bb')\otimes H_N
     + \sum_{C\ne H_N} Z_C\otimes C,
\]
where $C\in \dot{\tilde{\mathbf{B}}}$ runs over the canonical basis elements 
other than $H_N$, and $\iota:\dot{\mathbf{U}}1_{\zeta'}\rightarrow \dot{\tilde{\mathbf{U}}}1_{\tilde{\zeta}'}$ is an injective map determined by $\tilde{\zeta}$.
\end{theorem}

The map $\iota$ preserves the multiplication and sends $\dot{\mathbf{B}}\cap \dot{\mathbf{U}}1_{\zeta'}$ to $\dot{\tilde{\mathbf{B}}}\cap \dot{\tilde{\mathbf{U}}}1_{\tilde{\zeta}'}$. If $\dot{b}\dot{b}'=\sum_{\dot{b}''\in \dot{\mathbf{B}}}\dot{m}_{\dot{b},\dot{b}'}^{\dot{b}''}\dot{b}''$, then the sum of all terms in $\dot{\Delta}_{\tilde\zeta',0}(B_N)$ with the second factor $H_N$ is precisely 
\[
 \sum_{\dot{b}''\in \dot{\mathbf{B}}} \dot{m}_{\dot{b},\dot{b}'}^{\dot{b}''}\,\iota(\dot{b}'')\otimes H_N .
\]
In particular, all coefficients in this product $\dot{b}\dot{b}'$ occur in the coproduct of $\dot{\Delta}_{\tilde\zeta',0}(B_N)$ of the single canonical basis element $B_N$. For example, all nine terms of the multiplication in type $A_2$ in
\S\ref{subsec:a2-a4} occur in a single comultiplication in type $A_4$.
The converse fails already in type $A_1$: $v$ occurs as a structure
constant of comultiplication, whereas every structure constant of
multiplication is fixed by $v\mapsto v^{-1}$ (\S\ref{subsec:converse}).

The idea can already be seen in the negative part $\mathbf U^-$.
For canonical basis elements $b,b'$ of $\mathbf U^-$, we insert
$\theta_N=\prod_{i\in I}\theta_{i'}^{(N)}$ between them. In the
coproduct $r(b\theta_N b')\in
\tilde{\mathbf U}^-\otimes\tilde{\mathbf U}^-$, consider the terms
whose second factor has degree $|\theta_N|$. Only the new vertices
occur in this degree, so all the divided powers in $\theta_N$ go to
the second factor and all of $b,b'$ go to the first. This gives
$bb'\otimes\theta_N$ up to a power of $v$.

In the modified quantum group, we insert both positive and negative
divided powers at the new vertices. Let $D_N=\prod_{i\in I}E_{i'}^{(N)}F_{i'}^{(N)}$. Then we insert $D_N1_{\tilde\zeta}$ between the the lifts of the canonical basis elements $\dot b\in \dot{\mathbf{B}}\cap \dot{\mathbf{U}}1_{\zeta}$ and $\dot b'\in \dot{\mathbf{B}}\cap \dot{\mathbf{U}}1_{\zeta'}$, where $\tilde\zeta$ is a suitable lift of $\zeta$. In the coproduct of the resulting element $B_N$, we select the terms whose second factor is the canonical basis element $H_N=D_N1_0$. Theorem~\ref{thm:main}
shows that their sum is $\iota(\dot b\dot b')\otimes H_N$ for
every $N$. Proposition~\ref{prop:product} shows that the element
$B_N$  is a canonical
basis element of the larger modified quantum group for sufficiently large $N$.

Since the larger Cartan datum is again symmetric, our positivity
theorem for comultiplication \cite[Theorem~6.3]{Fang-He-2026} applies
to $B_N$. Then Theorem~\ref{thm:intro-comparison} 
gives the positivity for multiplication in the original modified quantum
group. Passing from the simply connected root datum to an arbitrary
root datum as in \cite[\S~6.6]{Fang-He-2026}, we obtain the
multiplication assertion of Lusztig's conjecture.

\begin{theorem}[Multiplication positivity]
\label{thm:intro-multiplication-positivity}
For any root datum associated with a symmetric Cartan datum, the
structure constants of multiplication in the modified quantum group $\dot{\mathbf{U}}$ with respect to the canonical basis $\dot{\mathbf{B}}$ belong to $\mathbb{N}[v,v^{-1}]$.
\end{theorem}

The multiplication conjecture was previously known for simply laced
finite type and for affine type $A$. For finite type $A$ and affine type $A$, it follows from a series work by Beilinson-Lusztig-MacPherson \cite{BLM}, Ginzburg-Vasserot \cite{Ginsburg-Vasserot-1993}, Lusztig \cite{Lusztig-1999,Lusztig-2000}, Schiffmann-Vasserot \cite{SV}, McGerty \cite{McGerty} and Fu-Shoji \cite{Fu-Shoji}. Fu-Shoji \cite{Fu-Shoji} identified the structure constants of multiplication in the modified quantum group $\dot{\mathbf{U}}(\widehat{\mathfrak{sl}}_n)$ with certain structure constants of multiplication in the positive part of a larger quantum group $\mathbf{U}(\widehat{\mathfrak{sl}}_m)$ with $n<m$. In any simply laced finite types, it follows from Webster's categorification \cite{W}. For arbitrary symmetric Cartan data, our earlier work established multiplication positivity when one factor belongs to a standard or
opposite parabolic subalgebra with Levi factor of finite type \cite[Theorem~6.6]{Fang-He-2026}. The present theorem removes this restriction on the factors.

{\bf Acknowledgements} JF is partially supported by the New Cornerstone Science Foundation through the New Cornerstone Investigator Program awarded to XH and the National Natural Science Foundation of China (Grant No. 12471030). XH is partially supported by the Hong Kong RGC grant 14300023 and by the New Cornerstone Science Foundation through the New Cornerstone Investigator Program.

\section{Preliminaries}\label{sec:framing}

\subsection{Quantum groups}
Let $I$ be a finite set, and let $A=(a_{ij})_{i,j\in I}$ be a
symmetrizable generalized Cartan matrix. We fix positive integers $d_i$ for any $i\in I$ such that $d_i a_{ij}=d_j a_{ji}$, and write $i\cdot j=d_i a_{ij}$ for any $i,j\in I$. Then $(I,\cdot)$ is a Cartan datum (see \cite[\S 1.1.1]{Lusztig-1993}). We call the Cartan datum $(I,\cdot)$ symmetric if $d_i=1$ for any $i\in I$.

Let $(Y,X,\langle\,,\,\rangle,\ldots)$ be the simply connected root datum of type $(I,\cdot)$, where $Y=\mathbb{Z}[I], X=\mathrm{Hom}(Y,\mathbb Z)$ and $\langle\,,\,\rangle:Y\times X\rightarrow \mathbb{Z}$ is the natural evaluation (see \cite[\S 2.2.2]{Lusztig-1993}). We denote by $i_Y$ and $i_X$ the simple coroot and the simple root for any $i\in I$. Then $\langle i_Y,j_X\rangle=a_{ij}$ for any $i,j\in I$. 

For any $\nu=\sum_{i\in I}\nu_i i\in\mathbb{Z}[I]$, any $n\in \mathbb{Z}$ and $l\leqslant m$ in $\mathbb{N}$, we denote
\begin{align*}
&\nu_Y=\sum_{i\in I}\nu_ii_Y\in Y,\ \nu_X=\sum_{i\in I}\nu_i i_X\in X,\ \mathrm{tr}\,\nu=\sum_{i\in I}\nu_i\in \mathbb{Z};\\
&[n]=\frac{v^n-v^{-n}}{v-v^{-1}},\ [m]!=\prod_{k=1}^m[k],\ \begin{bmatrix}
m\\l
\end{bmatrix}=\frac{[m]!}{[l]![m-l]!}\in \mathbb{Q}(v).
\end{align*}
For any $f(v)\in \mathbb{Q}(v)$ and $i\in I$, we denote $f_i(v)=f(v^{d_i})$.

The quantum group $\mathbf{U}$ associated with $(Y,X,\langle\,,\,\rangle,\ldots)$ is the $\mathbb{Q}(v)$-algebra with $1$ generated by $E_i,F_i$ and $K_{\mu}$ for any $i\in I$ and $\mu\in Y$, subject to the relations in \cite[\S 3.1.1]{Lusztig-1993}. The algebra $\mathbf{U}=\bigoplus_{\nu\in \mathbb{Z}[I]}\mathbf{U}(\nu)$ is $\mathbb{Z}[I]$-graded such that $E_i\in \mathbf{U}(i),F_i\in \mathbf{U}(-i)$ and $K_\mu\in \mathbf{U}(0)$ for any $i\in I$ and $\mu\in Y$. The tensor product $\mathbf{U}\otimes \mathbf{U}$ over $\mathbb{Q}(v)$ is an algebra with the natural multiplication $(u_1\otimes u_2)(u'_1\otimes u'_2)=u_1u'_1\otimes u_2u'_2$ for any  $u_1,u_2,u'_1,u'_2\in \mathbf{U}$.
The comultiplication of $\mathbf{U}$ is the algebra homomorphism $\Delta:\mathbf{U}\rightarrow \mathbf{U}\otimes \mathbf{U}$ defined by $$E_i\mapsto E_i\otimes 1+K_{d_ii_Y}\otimes E_i, F_i\mapsto 1\otimes F_i+F_i\otimes K_{-d_ii_Y}, K_\mu\mapsto K_\mu\otimes K_\mu\ \textrm{for any}\ i\in I,\mu\in Y.$$
Let $\mathbf{U}^+$ and $\mathbf{U}^-$ be the subalgebras of $\mathbf{U}$ generated by $E_i$ and $F_i$ for any $i\in I$ respectively. Then there is a natural algebra isomorphism $\mathbf{U}^+\rightarrow \mathbf{U}^-$ defined by $E_i\mapsto F_i$ for any $i\in I$. To study $\mathbf{U}^+$ and $\mathbf{U}^-$ uniformly, we follow \cite[Chapter 1]{Lusztig-1993} to introduce the algebra $\mathbf{f}$ associated with $(I,\cdot)$. 

The algebra $\mathbf{f}$ is the $\mathbb{Q}(v)$-algebra with $1$ generated by $\theta_i$ for any $i\in I$, subject to the quantum Serre relations in \cite[\S 1.4]{Lusztig-1993}. Then there are algebra isomorphisms $\mathbf{f}\rightarrow \mathbf{U}^+,x\mapsto x^+$ and $\mathbf{f}\rightarrow \mathbf{U}^-,x\mapsto x^-$ defined by $\theta_i^+=E_i$ and $\theta_i^-=F_i$ for any $i\in I$ respectively. The algebra $\mathbf{f}=\bigoplus_{\nu\in \mathbb{N}[I]}\mathbf{f}_{\nu}$ is $\mathbb{N}[I]$-graded such that $\theta_i\in \mathbf{f}_i$ for any $i\in I$. For any homogeneous $x\in \mathbf{f}_{\nu}$, we denote $|x|=\nu$. The tensor product $\mathbf{f}\otimes \mathbf{f}$ over $\mathbb{Q}(v)$ is an algebra with the twisted multiplication 
$(x_1\otimes x_2)(y_1\otimes y_2)=v^{|x_2|\cdot |y_1|}x_1y_1\otimes x_2y_2$ for any homogeneous $x_1,x_2,y_1,y_2\in \mathbf{f}$. The comultiplication of $\mathbf{f}$ is the algebra homomorphism $r:\mathbf{f}\rightarrow \mathbf{f}\otimes \mathbf{f}$ defined by $$\theta_i\mapsto \theta_i\otimes 1+1\otimes \theta_i\ \textrm{for any}\ i\in I.$$

Let $\dot{\mathbf{U}}$ be the modified quantum group associated with $(Y,X,\langle\,,\,\rangle,\ldots)$ defined in \cite[\S 23.1]{Lusztig-1993}. We have $\dot{\mathbf{U}}=\bigoplus_{\gamma,\zeta\in X}{_{\gamma}\mathbf{U}_{\zeta}}$, where
$${_{\gamma}\mathbf{U}_{\zeta}}=\mathbf{U}/(\sum_{\mu\in Y}(K_{\mu}-v^{\langle \mu,\gamma\rangle})\mathbf{U}+\sum_{\mu\in Y}\mathbf{U}(K_{\mu}-v^{\langle \mu,\zeta\rangle}))\ \textrm{for any}\ \gamma,\zeta\in X.$$
Let ${_{\gamma}\pi_{\zeta}}:\mathbf{U}\rightarrow {_{\gamma}\mathbf{U}_{\zeta}}$ be the natural projection for any $\gamma,\zeta\in X$, and let $1_{\zeta}={_{\zeta}\pi_{\zeta}}(1)$ for any $\zeta\in X$. The multiplication of $\mathbf{U}$ induces a well-defined multiplication on $\dot{\mathbf{U}}$ such that 
$${_{\gamma}\pi_{\zeta}}(u)\cdot {_{\gamma'}\pi_{\zeta'}}(u')=\delta_{\zeta,\gamma'}\cdot {_{\gamma}\pi_{\zeta'}}(uu')$$ 
for any $\gamma,\gamma',\zeta,\zeta'\in X,\nu,\nu'\in \mathbb{Z}[I]$ and $u,u'\in \mathbf{U}$ such that $\gamma-\zeta=\nu_X,\gamma'-\zeta'=\nu'_X$ and $u\in \mathbf{U}(\nu),u'\in \mathbf{U}(\nu')$. The comultiplication  of $\mathbf{U}$ induces a well-defined linear map ${_{\gamma',\gamma''}\dot{\Delta}_{\zeta',\zeta''}}:{_{\gamma'+\gamma''}\mathbf{U}_{\zeta'+\zeta''}}\rightarrow {_{\gamma'}\mathbf{U}_{\zeta'}}\otimes {_{\gamma''}\mathbf{U}_{\zeta''}}$ for any $\gamma',\gamma'',\zeta',\zeta''\in X$ such that 
$${_{\gamma',\gamma''}\dot{\Delta}_{\zeta',\zeta''}}({_{\gamma'+\gamma''}\pi_{\zeta'+\zeta''}}(u))=({_{\gamma'}\pi_{\zeta'}}\otimes {_{\gamma''}\pi_{\zeta''}})(\Delta(u))\ \textrm{for any}\ u\in \mathbf{U}.$$
The comultiplication of $\dot{\mathbf{U}}$ is the collection $\{{_{\gamma',\gamma''}\dot{\Delta}_{\zeta',\zeta''}};\gamma',\gamma'',\zeta',\zeta''\in X\}$.

Let $X^+=\{\lambda\in X;\langle i_Y,\lambda\rangle\in \mathbb{N}\ \textrm{for any}\ i\in I\}$ be the set of dominant weights. For any $\lambda\in X^+$, let ${^{\omega}\Lambda_{\lambda}}$ be the simple lowest weight module of $\mathbf{U}$ with the lowest weight $-\lambda$ and the lowest weight vector $\xi_{-\lambda}$, and let $\Lambda_{\lambda}$ be the simple highest weight module of $\mathbf{U}$ with the highest weight $\lambda$ and the highest weight vector $\eta_{\lambda}$ (see \cite[\S 3.4]{Lusztig-1993}). We regard them as $\dot{\mathbf{U}}$-modules via \cite[\S 23.1.4]{Lusztig-1993}.

\subsection{Canonical bases}

Let $\mathbf{B}$ be the canonical basis of the algebra $\mathbf{f}$ defined in \cite[Definition 14.4.6]{Lusztig-1993}. Under the isomorphisms $\mathbf{f}\rightarrow \mathbf{U}^+$ and $\mathbf{f}\rightarrow \mathbf{U}^-$, the canonical basis $\mathbf{B}$ of $\mathbf{f}$ gives the canonical bases $\{b^+;b\in \mathbf{B}\}$ of $\mathbf{U}^+$ and $\{b^-;b\in \mathbf{B}\}$ of $\mathbf{U}^-$ respectively. 

For any $\lambda\in X^+$, let $\mathbf{B}({^{\omega}\Lambda_{\lambda}})$ and $\mathbf{B}(\Lambda_{\lambda})$ be the canonical bases of modules ${^{\omega}\Lambda_{\lambda}}$ and $\Lambda_{\lambda}$ defined in \cite[Definition 14.4.12]{Lusztig-1993} respectively. By \cite[Theorem 14.4.11]{Lusztig-1993}, we have
$$\mathbf{B}({^{\omega}\Lambda_{\lambda}})=\{b^+\xi_{-\lambda};b\in \mathbf{B}\}\setminus \{0\},\qquad \mathbf{B}(\Lambda_{\lambda})=\{b^-\eta_{\lambda};b\in \mathbf{B}\}\setminus \{0\}.$$

For any $\lambda_1,\lambda_2\in X^+$, let $\mathbf{B}({^{\omega}\Lambda_{\lambda_1}}\otimes \Lambda_{\lambda_2})=\{b_1\diamondsuit b_2;b_1\in \mathbf{B}({^{\omega}\Lambda_{\lambda_1}}),b_2\in \mathbf{B}(\Lambda_{\lambda_2})\}$ be the canonical basis of the tensor product ${^{\omega}\Lambda_{\lambda_1}}\otimes \Lambda_{\lambda_2}$ defined in \cite[\S 24.3.4]{Lusztig-1993}. By \cite[Theorem 24.3.3]{Lusztig-1993}, for any $b_1\in \mathbf{B}({^{\omega}\Lambda_{\lambda_1}})$ and $b_2\in \mathbf{B}(\Lambda_{\lambda_2})$, the canonical basis element $b_1\diamondsuit b_2$ is the unique element of $\mathbb{Z}[v^{-1}][\mathbf{B}({^{\omega}\Lambda_{\lambda_1}})\otimes \mathbf{B}(\Lambda_{\lambda_2})]$ such that 
\begin{equation}\label{24.3.3}
\Psi(b_1\diamondsuit b_2)=b_1\diamondsuit b_2\ \textrm{and}\ b_1\diamondsuit b_2-b_1\otimes b_2\in v^{-1}\mathbb{Z}[v^{-1}][\mathbf{B}({^{\omega}\Lambda_{\lambda_1}})\otimes \mathbf{B}(\Lambda_{\lambda_2})],
\end{equation}
where $\Psi:{^{\omega}\Lambda_{\lambda_1}}\otimes \Lambda_{\lambda_2}\rightarrow {^{\omega}\Lambda_{\lambda_1}}\otimes \Lambda_{\lambda_2}$ is a $\mathbb{Q}$-linear involution defined in \cite[\S 24.3.2]{Lusztig-1993}.

Let $\dot{\mathbf{B}}=\{b_1\diamondsuit_{\zeta}b_2;\zeta\in X, b_1,b_2\in \mathbf{B}\}$ be the canonical basis of the modified quantum group $\dot{\mathbf{U}}$ defined in \cite[\S 25.2]{Lusztig-1993}. By \cite[Theorem 25.2.1]{Lusztig-1993}, for any $\zeta\in X$ and $b_1,b_2\in \mathbf{B}$, the canonical basis element $b_1\diamondsuit_{\zeta}b_2$ is the unique element of the integral form of $\dot{\mathbf{U}}1_{\zeta}$ such that for any $\lambda_1,\lambda_2\in X^+$, we have 
\begin{equation}\label{25.2.1}
(b_1\diamondsuit_{\zeta}b_2)(\xi_{-\lambda_1}\otimes \eta_{\lambda_2})=\begin{cases}
b_1^+\xi_{-\lambda_1}\diamondsuit b_2^-\eta_{\lambda_2}, \ &\textrm{if}\ \zeta=\lambda_2-\lambda_1\ \textrm{and}\ b_1^+\xi_{-\lambda_1}\neq 0, b_2^-\eta_{\lambda_2}\neq 0;\\
0, \ &\textrm{otherwise}.
\end{cases}    
\end{equation}

\section{Multiplication from comultiplication in quantum groups}\label{Structure constants}

In this section, we show that the structure constants of multiplication in the modified quantum group $\dot{\mathbf{U}}$ can be obtained from the structure constants of comultiplication of an explicit element in a larger modified quantum group $\dot{\tilde{\mathbf{U}}}$.

\subsection{Thickening data}
Let $\tilde{I}=I\sqcup I'$ be the finite set obtained from $I$ by adding a copy $I'=\{i';\ i\in I\}$ of $I$, and let $\tilde{A}=\begin{pmatrix}A&-\mathrm{Id}\\-\mathrm{Id}&2\mathrm{Id}\end{pmatrix}_{\tilde{I}\times \tilde{I}}$ be the symmetrizable generalized Cartan matrix extended from $A$ with $d_{i'}=d_i$ for any $i\in I$. The corresponding Cartan datum $(\tilde{I},\cdot)$ is defined to be the thickening Cartan datum of $(I,\cdot)$. Then the simply connected root datum $(\tilde{Y},\tilde{X},\langle\,,\,\rangle,\ldots)$ of type $(\tilde{I},\cdot)$ is defined to be the thickening root datum of $(Y,X,\langle\,,\,\rangle,\ldots)$.

When $(I,\cdot)$ is symmetric, these definitions coincide with those in \cite[\S 2.1]{Fang-He-2026}.

\subsection{Structure constants in negative parts}
\label{sec:negative}

In this subsection, we show that the structure constants of multiplication in the algebra $\mathbf{f}$ can be obtained from the structure constants of comultiplication in the larger algebra $\tilde{\mathbf{f}}$. This motivates the comparison in \S \ref{sec:Structure constants}.

Let $\tilde{\mathbf{f}}$ be the algebra associated with the thickening Cartan datum, and let $\tilde{\mathbf{B}}$ be its canonical basis. Then the imbedding $I\hookrightarrow \tilde{I}$ induces the subalgebra imbedding $\mathbf{f}\hookrightarrow \tilde{\mathbf{f}}$ and inclusion $\mathbf{B}\hookrightarrow \tilde{\mathbf{B}}$. 

For any $N\in \mathbb{N}$ and $i'\in I'$, we denote $\theta_{i'}^{(N)}=\theta_{i'}^{N}/([N]!)_{i'}$ and define $\theta_{N}=\prod_{i\in I}\theta_{i'}^{(N)}\in \tilde{\mathbf{f}}$. Since $i'\cdot j'=0$ and $\theta_{i'}\theta_{j'}=\theta_{j'}\theta_{i'}$ for any $i'\neq j'$ in $I'$, this element is independent of the order in the multiplication. By \cite[\S 14.5.3 \& Theorem 14.3.2(c)]{Lusztig-1993}, $\theta_{N}\in \tilde{\mathbf{B}}$ is a canonical basis element of $\tilde{\mathbf{f}}$.

For any $b,b'\in \mathbf{B}$ and $x\in \tilde{\mathbf{f}}$, we write $bb'=\sum_{b''\in \mathbf{B}}m^{b''}_{b,b'}b''$ and $r(x)=\sum_{b_1,b_2\in \tilde{\mathbf{B}}}c_{b_1,b_2}^xb_1\otimes b_2$.

\begin{proposition}\label{prop:half}
For any $b,b',b''\in \mathbf{B}$ and $N\in \mathbb{N}$, we have
$$m_{b,b'}^{b''}=v^{-|\theta_N|\cdot |b'|}c_{b'',\theta_N}^{b\theta_N b'}.$$
\end{proposition}
\begin{proof}
Consider $(\mathrm{Id}\otimes p)\,r:\tilde{\mathbf{f}}\rightarrow \tilde{\mathbf{f}}\otimes \tilde{\mathbf f}_{|\theta_N|}$, where $p:\tilde{\mathbf f}\rightarrow \tilde{\mathbf f}_{|\theta_N|}$ be the
projection onto the homogeneous component of degree $|\theta_N|$.
Since this degree is supported on the new vertices in $I'$, no generator
from $b$ or $b'$ can occur in the second factor, and every generator from $\theta_N$ must occur there. Thus the only
terms of $r(b)$, $r(\theta_N)$ and $r(b')$ that contribute are
$b\otimes 1$, $1\otimes\theta_N$ and $b'\otimes 1$ respectively.
We obtain
\[
 (\mathrm{Id}\otimes p)\,r(b\theta_N b')
   =(b\otimes 1)(1\otimes\theta_N)(b'\otimes 1)
   =v^{|\theta_N|\cdot|b'|}\,bb'\otimes\theta_N .
\]
Comparing the coefficients of $b''\otimes\theta_N$ on both sides
gives the assertion.
\end{proof}

For a symmetric Cartan datum, we will show that $b\theta_N b'$ is a canonical basis element of $\tilde{\mathbf{f}}$ for sufficiently large $N$ in Lemma \ref{lem:sep}.

\subsection{Lifts of canonical bases}\label{Imbedding of canonical bases}

Let $\tilde{\mathbf{U}}$ and $\dot{\tilde{\mathbf{U}}}$ be the quantum group and modified quantum group associated with the thickening root datum respectively. Then the imbedding $I\hookrightarrow \tilde{I}$ induces the subalgebra imbedding $\mathbf{U}\hookrightarrow \tilde{\mathbf{U}}$. Let $\tilde\zeta\in\tilde X$ be a lift of $\zeta\in X$. Since $\zeta$ and $\tilde{\zeta}$ agree on the original coweight lattice $Y$, by the triangular decomposition of modified quantum group in \cite[\S 23.2.1]{Lusztig-1993}, there is a well-defined injective linear map
$$\iota_{\zeta,\tilde\zeta}:
\dot{\mathbf U}1_\zeta\longrightarrow \dot{\tilde{\mathbf U}}1_{\tilde\zeta},\, x^+y^-1_\zeta\longmapsto x^+y^-1_{\tilde\zeta}\ \textrm{for any}\ x,y\in\mathbf{f}.$$
Let $\dot{\tilde{\mathbf{B}}}$ be the canonical basis of $\dot{\tilde{\mathbf{U}}}$.

\begin{lemma}\label{lift}
For any $b_1,b_2\in \mathbf{B}$, we have $\iota_{\zeta,\tilde{\zeta}}(b_1\diamondsuit_{\zeta} b_2)=b_1\diamondsuit_{\tilde{\zeta}} b_2$.
\end{lemma}

\begin{remark}
We call $b_1\diamondsuit_{\tilde{\zeta}} b_2\in \dot{\tilde{\mathbf{B}}}$ a {\it lift} of $b_1\diamondsuit_{\zeta}b_2\in \dot{\mathbf{B}}$.
\end{remark}

\begin{proof}
Choose $\tilde\lambda_1,\tilde\lambda_2\in\tilde X^+$ with
$\tilde\zeta=\tilde\lambda_2-\tilde\lambda_1$, and let
$\lambda_1,\lambda_2\in X^+$ be their restrictions to the original coweight lattice $Y$. So $\langle \mu,\lambda_t\rangle=\langle \mu,\tilde{\lambda}_t\rangle$ for any $\mu\in Y$ and $t=1,2$. Then $\zeta=\lambda_2-\lambda_1$. We denote the generators of the tensor products
\[
 c=\xi_{-\lambda_1}\otimes\eta_{\lambda_2}\in {^{\omega}\Lambda_{{\lambda}_1}}\otimes \Lambda_{{\lambda}_2},
 \qquad
 \tilde{c}=\xi_{-\tilde\lambda_1}\otimes\eta_{\tilde\lambda_2}\in {^{\omega}\Lambda_{\tilde{\lambda}_1}}\otimes \Lambda_{\tilde{\lambda}_2}.
\]
The $\mathbf U$-submodules generated by $\xi_{-\tilde\lambda_1}\in {^{\omega}\Lambda_{\tilde{\lambda}_1}}$
and $\eta_{\tilde\lambda_2}\in \Lambda_{\tilde{\lambda}_2}$ are naturally isomorphic to ${}^{\omega}\Lambda_{\lambda_1}$ and $\Lambda_{\lambda_2}$
respectively. We denote the corresponding inclusions by $\iota_1$ and $\iota_2$,
and denote $\iota=\iota_1\otimes \iota_2$. Then we have the commutative diagram
\[
\begin{array}{ccc}
 \dot{\mathbf U}1_\zeta
 & \xrightarrow{\ \iota_{\zeta,\tilde\zeta}\ }
 & \dot{\tilde{\mathbf U}}1_{\tilde\zeta} \\[4pt]
 {\scriptstyle u\mapsto uc}\,\big\downarrow
 &&
 \big\downarrow\,{\scriptstyle \tilde u\mapsto\tilde u\tilde c} \\[4pt]
 {}^{\omega}\Lambda_{\lambda_1}\otimes\Lambda_{\lambda_2}
 & \xrightarrow{\ \iota\ }
 & {}^{\omega}\Lambda_{\tilde\lambda_1}
      \otimes\Lambda_{\tilde\lambda_2}.
\end{array}
\]
Indeed, since $\iota$ is $\mathbf U$-linear and $\iota(c)=\tilde c$, both compositions give $u\tilde c$ on each $u1_\zeta\in \dot{\mathbf U}1_\zeta$ for $u\in \mathbf{U}$.

Let $\Psi$ and $\tilde\Psi$ denote the involutions on the source and target of $\iota$ respectively defined in \cite[\S 24.3.2]{Lusztig-1993}. The map $\tilde\Psi\circ \iota\circ\Psi$ is $\mathbf U$-linear and
sends $c$ to $\tilde c$. Since $c$ generates the source of $\iota$, this map equals $\iota$, and so $\tilde\Psi\circ \iota=\iota\circ\Psi$.
Moreover, $\mathbf B\subseteq\tilde{\mathbf B}$, so $\iota_1$ and $\iota_2$ send canonical basis elements to canonical basis elements. For $p\in\mathbf B({}^{\omega}\Lambda_{\lambda_1})$ and
$q\in\mathbf B(\Lambda_{\lambda_2})$, we have
\begin{gather*}
 \tilde\Psi\bigl(\iota(p\diamondsuit q)\bigr)
 =\iota\bigl(\Psi(p\diamondsuit q)\bigr)
 =\iota(p\diamondsuit q), \\
 \iota(p\diamondsuit q)-\iota_1(p)\otimes \iota_2(q)
 =\iota(p\diamondsuit q-p\otimes q) \in v^{-1}\mathbb Z[v^{-1}]
 \bigl[\mathbf B({}^{\omega}\Lambda_{\tilde\lambda_1})
       \otimes\mathbf B(\Lambda_{\tilde\lambda_2})\bigr].
\end{gather*}
By \eqref{24.3.3}, these two conditions give $\iota(p\diamondsuit q)=\iota_1(p)\diamondsuit \iota_2(q)$. Applying the diagram to $b_1\diamondsuit_\zeta b_2$ and using
\eqref{25.2.1}, we obtain
\begin{equation}\label{characterizing condition}
 \iota_{\zeta,\tilde\zeta}(b_1\diamondsuit_\zeta b_2)\,\tilde c
 =\iota\bigl((b_1\diamondsuit_\zeta b_2)c\bigr)
 =(b_1\diamondsuit_{\tilde\zeta}b_2)\,\tilde c .
\end{equation}
If $b_1^+\xi_{-\lambda_1}=0$ or $b_2^-\eta_{\lambda_2}=0$,
the corresponding vector in the enlarged module is also zero, so both
sides vanish.
As $(\tilde\lambda_1,\tilde\lambda_2)$
varies over the dominant pairs with  $\tilde\zeta=\tilde\lambda_2-\tilde\lambda_1$, and
action on $\tilde c$ separates the elements of $\dot{\tilde{\mathbf U}}1_{\tilde\zeta}$ by \eqref{25.2.1}.
The lemma follows.
\end{proof}

\subsection{The \texorpdfstring{$\star$}{$\star$}-product}\label{Auxiliary elements}

Let $\dot b=b_1\diamondsuit_\zeta b_2,\dot b'=b'_1\diamondsuit_{\zeta'}b'_2\in \dot{\mathbf B}$. Then we denote $|\dot{b}|=|b_1|-|b_2|,|\dot{b}'|=|b'_1|-|b'_2|\in\mathbb Z[I]$. We may assume that $\zeta=\zeta'+|\dot b'|_X$, since otherwise $\dot b\dot b'=0$.

Let $\tilde\zeta,\tilde\zeta'\in \tilde{X}$ be the lifts of $\zeta,\zeta'\in X$ determined by $\langle i'_{\tilde Y},\tilde\zeta\rangle=0,\langle i'_{\tilde Y},\tilde\zeta'\rangle=|\dot b'|_i$ for any $i\in I$ respectively. With this choice, $\tilde\zeta'+|\dot b'|_{\tilde X}=\tilde\zeta$,
so the lift $b'_1\diamondsuit_{\tilde{\zeta}'}b'_2$ of $\dot b'$ lies in $1_{\tilde{\zeta}}\dot{\tilde{\mathbf{U}}}1_{\tilde{\zeta'}}$. We denote $\tilde{\gamma}=\tilde{\zeta}+|\dot{b}|_{\tilde{X}}\in \tilde{X}$. Then the lift $b_1\diamondsuit_{\tilde{\zeta}}b_2$ of $\dot{b}$ lies in $1_{\tilde{\gamma}}\dot{\tilde{\mathbf{U}}}1_{\tilde{\zeta'}}$.

Fix $N\in\mathbb N$, and let $\theta_{N}=\prod_{i\in I}\theta_{i'}^{(N)}\in \tilde{\mathbf{f}}$ be as in
\S\ref{sec:negative}. We define the \emph{$\star$-product} of
$\dot b$ and $\dot b'$ by inserting
$\theta_N\diamondsuit_{\tilde\zeta}\theta_N$ between their lifts:
\begin{equation}\label{star}
 \dot b \star \dot b'
 :=
 (b_1\diamondsuit_{\tilde\zeta}b_2)\,
 (\theta_N\diamondsuit_{\tilde\zeta}\theta_N)\,
 (b'_1\diamondsuit_{\tilde\zeta'}b'_2)
 \in 1_{\tilde{\gamma}}\dot{\tilde{\mathbf U}}1_{\tilde\zeta'}.
\end{equation}
We suppress $N$ in the notation.

\subsection{Structure constants}\label{sec:Structure constants}

We compare the expansion of $\dot b\dot b'$ in the canonical basis with the term in the comultiplication of $\dot b \star \dot b'$ whose second factor is $\theta_N\diamondsuit_0\theta_N$.

For the multiplication in $\dot{\mathbf{U}}$, we write
\[
 \dot b\dot b'
   =\sum
     \dot m_{\dot b,\dot b'}^{\dot b''}\,\dot b'',
\]
where $\dot b''\in\dot{\mathbf B}\cap1_{\tilde{\gamma}}\dot{\mathbf U}1_{\zeta'}$. For the comultiplication in $\dot{\tilde{\mathbf{U}}}$, we write
\[
 _{\tilde{\gamma},0}\dot\Delta_{\tilde\zeta',0}(\dot b \star \dot b')
   =\sum \dot{\tilde c}_{\dot{b}_1,\dot{b}_2}^{\,\dot b \star \dot b'}\,\dot{b}_1\otimes \dot{b}_2,
\]
where $\dot{b}_1\in\dot{\tilde{\mathbf B}}\cap
1_{\tilde{\gamma}}\dot{\tilde{\mathbf U}}1_{\tilde\zeta'}$ and
$\dot{b}_2\in \dot{\tilde{\mathbf B}}\cap 1_0\dot{\tilde{\mathbf U}}1_0$.
Both summations are finite.

We denote $H_N=\theta_N\diamondsuit_0\theta_N$, and let 
$$H_N^*:\dot{\tilde{\mathbf{U}}}:\rightarrow \mathbb{Q}(v), x\mapsto \delta_{x,H_N}\ \textrm{for any}\ x\in \tilde{\mathbf{B}}$$
be the corresponding dual canonical basis element.

\begin{lemma}\label{unique term}
For $u_1,u_2\in\mathbf U$ and $z\in\dot{\tilde{\mathbf U}}$, we have
\[
 H_N^*(u_1zu_2)=\mathbf e(u_1)\,H_N^*(z)\,\mathbf e(u_2),
\]
where $\mathbf e$ is the counit of $\mathbf U$.
\end{lemma}

\begin{proof}
Let $\tilde\lambda=N\sum_{i\in I}\omega_{i'}$, where $\omega_{i'},i\in I$ are fundamental weights. Let $c=\xi_{-\tilde\lambda}\otimes\eta_{\tilde\lambda}\in {^{\omega}\Lambda_{\tilde\lambda}\otimes \Lambda_{\lambda}}$. Since $\theta_N^+\xi_{-\tilde\lambda}\neq 0$ and $\theta_N^-\eta_{\tilde\lambda}\neq 0$, \eqref{25.2.1} gives
\[
 H_Nc=(\theta_N^+\xi_{-\tilde\lambda})
       \diamondsuit(\theta_N^-\eta_{\tilde\lambda}).
\]
Moreover, the nonzero images of the canonical basis elements in $\dot{\tilde{\mathbf{B}}}$ under the actions on $c$, are distinct canonical basis
elements in $\mathbf{B}({^{\omega}\Lambda_{\tilde\lambda}\otimes \Lambda_{\tilde{\lambda}}})$. Thus, for $z\in\dot{\tilde{\mathbf U}}$, we have
\[
 H_N^*(z)=(H_Nc)^*(zc),
\]
where $(H_Nc)^*:{^{\omega}\Lambda_{\tilde\lambda}\otimes \Lambda_{\tilde{\lambda}}}\rightarrow \mathbb{Q}(v)$ is the dual canonical basis element corresponding to $H_Nc$.

Since $\tilde\lambda$ vanishes on the original coweight lattice $Y$, the lines spanned by $\xi_{-\tilde\lambda}\in {^{\omega}\Lambda_{\tilde\lambda}}$ and
$\eta_{\tilde\lambda}\in \Lambda_{\tilde{\lambda}}$ are trivial $\mathbf U$-modules. Hence $uc=\mathbf e(u)c$ for any $u\in\mathbf U$, and then
\begin{equation}\label{right action}
H_N^*(zu)=(H_Nc)^*\bigl((zu)c\bigr)=(H_Nc)^*\bigl(z(uc)\bigr)=(H_Nc)^*(\mathbf{e}(u)zc)=\mathbf e(u)H_N^*(z).    
\end{equation}

For the left action, we use the automorphism $\sigma:\dot{\tilde{\mathbf{U}}}\rightarrow \dot{\tilde{\mathbf{U}}}^{\mathrm{opp}}$.
Since the new vertices $i'\in I'$ are pairwise orthogonal, applying
\cite[Proposition~25.3.2]{Lusztig-1993} with $a=b=N$ and $n=2N$
at each vertex $i'$ gives
\[
 H_N
 =\Bigl(\prod_{i\in I}E_{i'}^{(N)}F_{i'}^{(N)}\Bigr)1_0
 =\Bigl(\prod_{i\in I}F_{i'}^{(N)}E_{i'}^{(N)}\Bigr)1_0 .
\]
Thus $\sigma(H_N)=H_N$. Since $\sigma(\dot{\tilde{\mathbf{B}}})=\dot{\tilde{\mathbf{B}}}$ by \cite[Theorem 4.3.2]{Kashiwara-1994}, we have
$k\circ\sigma=k$.
Applying the identity \eqref{right action} and using
$\mathbf e\circ\sigma=\mathbf e$, we obtain
\begin{equation}\label{left action}
H_N^*(uz)=H_N^*\bigl(\sigma(z)\sigma(u)\bigr)
 =\mathbf e\bigl(\sigma(u)\bigr)H_N^*\bigl(\sigma(z)\bigr)
 =\mathbf e(u)H_N^*(z).   
\end{equation}
The two identities \eqref{right action} and \eqref{left action} give the assertion.
\end{proof}

\begin{theorem}\label{thm:main}
For any $N\in\mathbb N$, we have
\[
 (\mathrm{Id}\otimes H_N^*)\,{_{\tilde{\gamma},0}
 \dot\Delta_{\tilde\zeta',0}}(\dot b \star \dot b')
 =\iota_{\zeta',\tilde\zeta'}(\dot b\dot b').
\]
In particular, for any $\dot b''\in\dot{\mathbf B}\cap\dot{\mathbf U}1_{\zeta'}$, we have $\dot m_{\dot b,\dot b'}^{\dot b''}=\dot{\tilde c}_{\iota_{\zeta',\tilde\zeta'}(\dot b''),\,H_N}^{\,\dot b \star \dot b'}$.
\end{theorem}
\begin{proof}
For $z\in\tilde{\mathbf U}$, set $\widehat k(z)=H_N^*(z1_0)=(H_Nc)^*(zc)$, with $c,(H_Nc)^*$ as in the proof of Lemma~\ref{unique term}. For $u,u'\in\mathbf U$, we have $u'c=\mathbf e(u')c$, while Lemma~\ref{unique term} gives $H_N^*(uz1_0)=\mathbf e(u)H_N^*(z1_0)$. Thus
\[
 \widehat k(uzu')=\mathbf e(u)\,\widehat k(z)\,\mathbf e(u').
\]
Since $\Delta$ is an algebra homomorphism and
$(\mathrm{Id}\otimes\mathbf e)\Delta(u)=u$ for any $u\in\mathbf U$, suppose that $\Delta(z)=\sum z_1\otimes z_2$, we have
\begin{equation}\label{Hopf1}
 (\mathrm{Id}\otimes\widehat k)\Delta(uzu')
 =\sum \widehat{k}(z_2)uz_1u'=u\,\bigl((\mathrm{Id}\otimes\widehat k)\Delta(z)\bigr)\,u'.
\end{equation}
Let 
\[
 D_N=\theta_N^+\theta_N^-
     =\prod_{i\in I}\bigl(E_{i'}^{(N)}F_{i'}^{(N)}\bigr)
     \in\tilde{\mathbf U}.
\]
Since the new vertices $i'\in I'$ are pairwise orthogonal and
$\langle i'_{\tilde Y},\tilde\zeta\rangle=0$ for all $i\in I$, by
\cite[Proposition~25.3.2]{Lusztig-1993}, we have
$D_N1_{\tilde\zeta}=\theta_N\diamondsuit_{\tilde\zeta}\theta_N$ and $D_N1_0=H_N$. 

We calculate $(\mathrm{Id}\otimes\widehat k)\Delta(D_N)$. For $\mathbf r=(r_i)_{i\in I}\in\mathbb N^I$, set $\theta_{\mathbf r}=\prod_{i\in I}\theta_{i'}^{(r_i)}$.
Applying the divided-power comultiplication formula (see \cite[\S3.1.5]{Lusztig-1993}), and multiplying the second factor on the right by $1_0$, we obtain
\[
 \Delta(D_N)(1\otimes 1_0)
 \in \sum_{0\leqslant r_i,s_i\leqslant N}
 \tilde{\mathbf U}\otimes
 \mathbb Q(v)\,\theta_{\mathbf r}^+\theta_{\mathbf s}^-1_0 .
\]
Since the new vertices are pairwise orthogonal,
$\theta_{\mathbf r}^+\theta_{\mathbf s}^-$ changes the
$i'$-weight by $2(r_i-s_i)$. Thus $1_0\theta_{\mathbf r}^+\theta_{\mathbf s}^-1_0=0$ unless $\mathbf r=\mathbf s$.
Only terms with $\mathbf r=\mathbf s$ can therefore contribute to
the coefficient of $H_N$. In that case, \cite[Proposition~25.3.2]{Lusztig-1993} gives $\theta_{\mathbf r}^+\theta_{\mathbf r}^-1_0
 =\theta_{\mathbf r}\diamondsuit_0\theta_{\mathbf r}$.
Consequently,
\[
 H_N^*(\theta_{\mathbf r}^+\theta_{\mathbf s}^-1_0)
 =\begin{cases}
     1, & \text{if } r_i=s_i=N \text{ for any } i\in I;\\
     0, & \text{otherwise}.
   \end{cases}
\]
The corresponding summand in
$\Delta(D_N)(1\otimes 1_0)$ is
$\bigl(\prod_{i\in I}K_{d_i i'_{\tilde Y}}^N\bigr)\otimes H_N$ with coefficient one. Hence
\[
 (\mathrm{Id}\otimes\widehat k)\Delta(D_N)
 =\prod_{i\in I}K_{d_i i'_{\tilde Y}}^N .
\]

Choose homogeneous $u,u'\in\mathbf U$ such that
$\dot b=u1_\zeta$ and $\dot b'=u'1_{\zeta'}$. Then
$\dot b \star \dot b'=uD_Nu'1_{\tilde\zeta'}$.
Since $u'1_{\tilde\zeta'}=1_{\tilde\zeta}u'1_{\tilde\zeta'}$ and
$\langle i'_{\tilde Y},\tilde\zeta\rangle=0$ for all $i\in I$, we have
\[
 \Bigl(\prod_{i\in I}K_{d_i i'_{\tilde Y}}^N\Bigr)
 u'1_{\tilde\zeta'}=u'1_{\tilde\zeta'}.
\]
Applying \eqref{Hopf1}, we obtain
\[
 \begin{aligned}
 (\mathrm{Id}\otimes H_N^*)\,
 {_{\tilde{\gamma},0}\dot\Delta_{\tilde\zeta',0}}(\dot b \star \dot b')
  =\bigl((\mathrm{Id}\otimes\widehat k)\Delta(uD_Nu')\bigr)
                   1_{\tilde\zeta'}
 =u\Bigl(\prod_{i\in I}K_{d_i i'_{\tilde Y}}^N\Bigr)
      u'1_{\tilde\zeta'}
 =uu'1_{\tilde\zeta'}
  =\iota_{\zeta',\tilde\zeta'}(\dot b\dot b').
 \end{aligned}
\]
Since $\iota_{\zeta',\tilde\zeta'}(\dot{\mathbf{B}})\subset \dot{\tilde{\mathbf{B}}}$ by Lemma \ref{lift}, comparison of coefficients gives the second assertion.
\end{proof}

\section{Canonicity and multiplication positivity}\label{sec:main}

Throughout this section, we assume that the Cartan datum $(I,\cdot)$
is symmetric. The thickened Cartan datum $(\tilde I,\cdot)$ is then
symmetric as well. We retain the notation of  \S\ref{Structure constants}.

Theorem \ref{thm:main} expresses the coefficients of the product $\dot b\dot b'$ in terms of the coproduct of $\dot b \star \dot b'$
for any $N\in \mathbb{N}$. We now show that $\dot b \star \dot b'$ is itself a
canonical basis element for all sufficiently large $N$.

\subsection{Canonicity in the negative part}

We first consider the element $b\theta_N b'\in \tilde{\mathbf f}$.
The following lemma gives the canonical basis elements used to describe
$\dot b \star \dot b'$ in Proposition \ref{prop:product}.

\begin{lemma}\label{lem:sep}
For any $b,b'\in\mathbf B$, we have $b\theta_N b'\in\tilde{\mathbf B}$
for all sufficiently large $N$.
\end{lemma}
\begin{proof}
The canonical basis elements $b,\theta_N,b'\in \tilde{\mathbf{B}}$ are integral and bar-invariant in the sense of the first two conditions in \cite[Theorem 14.2.3]{Lusztig-1993}, so is their product $b\theta_N b'\in \tilde{\mathbf{f}}$. We compute its norm with respect to the bilinear form on $\tilde{\mathbf f}$ defined in \cite[Proposition~1.2.3]{Lusztig-1993}, which restricts to the bilinear form on $\mathbf f$. Property (c) in \cite[Proposition~1.2.3]{Lusztig-1993} gives 
$$(b\theta_N b',b\theta_N b')=(b\otimes\theta_N\otimes b',(\mathrm{Id}\otimes r)r(b\theta_N b')).$$
Then the homogeneity of the bilinear form forces all generators at the new vertices in $I'$ into the middle factor and all generators at
the original vertices in $I$ into the two outer factors. More precisely, we write $r(b)=\sum x_1\otimes x_2$ and $r(b')=\sum y_1\otimes y_2$ with homogeneous factors. Then  
$$(b\theta_N b',b\theta_N b')=\sum v^{|x_2|\cdot|y_1|+(|x_2|+|y_1|)\cdot|\theta_N|}(b,x_1y_1)(\theta_N,\theta_N)(b',x_2y_2),$$
where the exponent comes from the twisted multiplication in the tensor product $\tilde{\mathbf{f}}^{\otimes 3}$. If $(b,x_1y_1)\neq 0$, then $|x_1|+|y_1|=|b|=|x_1|+|x_2|$, so
$|y_1|=|x_2|$ and
\[
 (b\theta_N b',b\theta_N b')
 =(\theta_N,\theta_N)
 \sum_{\nu\in\mathbb N[I]}
 v^{\nu\cdot\nu-2N\operatorname{tr}\nu}
 \sum_{|x_2|=|y_1|=\nu}
 (b,x_1y_1)(b',x_2y_2).
\]

The summand with $\nu=0$ is $(b,b)(b',b')$. For each $\nu\ne0$,
the inner summation is independent of $N$, whereas
$\nu\cdot\nu-2N\operatorname{tr}\nu\rightarrow -\infty$ when $N\rightarrow \infty$.
Since only finitely many $\nu$ occur, the sum over $\nu\ne0$ belongs to $v^{-1}\mathbb Q[[v^{-1}]]$ for all
sufficiently large $N$. The norms $(b,b),(b',b')$ and $(\theta_N,\theta_N)$
belong to $1+v^{-1}\mathbb Q[[v^{-1}]]$ by
\cite[Theorem~14.2.3]{Lusztig-1993}. Hence so does
$(b\theta_N b',b\theta_N b')$.

By almost orthonormality, an integral bar-invariant element
with such a norm belongs to $\pm\tilde{\mathbf B}$ by
\cite[Lemma~14.2.2 and Theorem~14.2.3]{Lusztig-1993}.
Since the thickening Cartan datum is symmetric, by \cite[Theorem~14.4.13(b)]{Lusztig-1993}, multiplication in
$\tilde{\mathbf f}$ has coefficients in $\mathbb N[v,v^{-1}]$ with
respect to $\tilde{\mathbf B}$. As $b,\theta_N,b'\in \tilde{\mathbf B}$, the negative sign is excluded.
Thus $b\theta_N b'\in\tilde{\mathbf B}$.
\end{proof}

\subsection{Canonicity in the modified quantum group}

Let $\dot b=b_1\diamondsuit_\zeta b_2$ and
$\dot b'=b'_1\diamondsuit_{\zeta'}b'_2$ be as in \S \ref{Auxiliary elements}.
By Lemma \ref{lem:sep}, the elements $b_1\theta_N b'_1$ and
$b_2\theta_N b'_2$ belong to $\tilde{\mathbf B}$ for all sufficiently
large $N$. We now show that their ``$\diamondsuit_{\tilde\zeta'}$-product''
is precisely $\dot b \star \dot b'$.

\begin{proposition}\label{prop:product}
For all sufficiently large $N$, we have
\[
 \dot b \star \dot b'
   =(b_1\theta_N b'_1)\diamondsuit_{\tilde\zeta'}
     (b_2\theta_N b'_2).
\]
In particular, $\dot b \star \dot b'$ is a canonical basis element of $\dot{\tilde{\mathbf U}}$.
\end{proposition}

\begin{proof}
Set $r_i=\max(|\dot b'|_i,0)$ for $i\in I$.
Choose $\lambda'_1,\lambda'_2\in X^+$ such that
$\zeta'=\lambda'_2-\lambda'_1$ and
$\langle i_Y,\lambda'_1\rangle\geq |b_1|_i+|b'_1|_i,\langle i_Y,\lambda'_2\rangle\geq |b_2|_i+|b'_2|_i$
for any $i\in I$. Take their lifts
$\tilde\lambda'_1,\tilde\lambda'_2\in\tilde X$ such that
\begin{equation}\label{lift lambda'_1,lambda'_2}
 \langle i'_{\tilde Y},\tilde\lambda'_1\rangle=N+r_i-|\dot b'|_i,\,
 \langle i'_{\tilde Y},\tilde\lambda'_2\rangle=N+r_i\geqslant N\ \textrm{for any}\ i\in I.
\end{equation}
The lifts $\tilde\lambda'_1,\tilde\lambda'_2\in\tilde X^+$ are dominant and satisfy $\tilde\zeta'=\tilde{\lambda}'_2-\tilde{\lambda}'_1$.

We simply denote $\xi=\xi_{-\tilde\lambda'_1},\eta=\eta_{\tilde\lambda'_2}$ and $c=\xi\otimes\eta\in {^{\omega}\Lambda_{\tilde{\lambda}'_1}}\otimes \Lambda_{\tilde{\lambda}'_2}$. Since each factor in $\dot{b}\star \dot{b}'=(b_1\diamondsuit_{\tilde{\zeta}} b_2)(\theta_N\diamondsuit_{\tilde{\zeta}} \theta_N)(b'_1\diamondsuit_{\tilde{\zeta}'} b'_2)$ is bar-invariant with respect to the bar-involution on $\dot{\tilde{\mathbf{U}}}$, $\dot{b}\star \dot{b}'$ is bar-invariant. Let $\tilde{\Psi}:{^{\omega}\Lambda_{\tilde{\lambda}'_1}}\otimes \Lambda_{\tilde{\lambda}'_2}\rightarrow {^{\omega}\Lambda_{\tilde{\lambda}'_1}}\otimes \Lambda_{\tilde{\lambda}'_2}$ be the involution defined in \cite[\S 24.3.2]{Lusztig-1993}. Then we have 
\begin{equation}\label{Psi-invariant}
\tilde{\Psi}\bigl((\dot b \star \dot b')c\bigr)
 =\overline{\dot b \star \dot b'}\,\tilde{\Psi}(c)
 =(\dot b \star \dot b')c .    
\end{equation}
Let $\mathcal L=\mathbb Z[v^{-1}][\mathbf{B}({^{\omega}\Lambda_{\tilde\lambda'_1}})\otimes \mathbf{B}(\Lambda_{\tilde\lambda'_2})]$.
It suffices to prove that, for any sufficiently large $N$, we have
\begin{equation}\label{eq:star-leading-tensor}
 (\dot b\star \dot b')c
 \equiv (b_1\theta_N b'_1)^+\xi\otimes(b_2\theta_N b'_2)^-\eta
 \mod{v^{-1}\mathcal L}.
\end{equation}
Indeed, by Lemma \ref{lem:sep} and the choices of $\tilde{\lambda}'_1,\tilde{\lambda}'_2$, both $(b_1\theta_N b'_1)^+\xi$ and $(b_2\theta_N b'_2)^-\eta$ are canonical basis elements. Thus, if \eqref{eq:star-leading-tensor} holds, then \eqref{Psi-invariant}, \eqref{24.3.3} and \eqref{25.2.1} give
\[
 \begin{aligned}
 (\dot b\star \dot b')c
 &= (b_1\theta_N b'_1)^+\xi\diamondsuit(b_2\theta_N b'_2)^-\eta= \bigl((b_1\theta_N b'_1)\diamondsuit_{\tilde\zeta'}
                 (b_2\theta_N b'_2)\bigr)c .
 \end{aligned}
\]
The difference $\dot b \star \dot b'-(b_1\theta_N b'_1)\diamondsuit_{\tilde\zeta'}(b_2\theta_N b'_2)$ lies in the kernel of the map $\dot{\tilde{\mathbf U}}1_{\tilde\zeta'}\rightarrow {^{\omega}\Lambda_{\tilde{\lambda}'_1}}\otimes \Lambda_{\tilde{\lambda}'_2}$ given by actions on $c$. By
\cite[Proposition 23.3.6 \& Theorem 25.2.5]{Lusztig-1993}, this kernel is
$$\sum_{\tilde{i}\in\tilde I}
 \Bigl(\dot{\tilde{\mathbf U}}\,
       E_{\tilde{i}}^{(\langle \tilde{i}_{\tilde Y},\tilde\lambda'_1\rangle+1)}
      +\dot{\tilde{\mathbf U}}\,
       F_{\tilde{i}}^{(\langle \tilde{i}_{\tilde Y},\tilde\lambda'_2\rangle+1)}\Bigr)
 1_{\tilde\zeta'}$$ 
spanned by $\dot{\tilde{\mathbf{B}}}$. If there exists $x_1\diamondsuit_{\tilde\zeta'}x_2\in \dot{\tilde{\mathbf{B}}}$ occurring in the linear combination of $\dot{\tilde{\mathbf{B}}}$ with nonzero coefficient for the difference $\dot b \star \dot b'-(b_1\theta_N b'_1)\diamondsuit_{\tilde\zeta'}(b_2\theta_N b'_2)$, then 
\[
 |x_a|_{\tilde{i}}\leqslant |b_a|_{\tilde{i}}+|\theta_N|_{\tilde{i}}+|b'_a|_{\tilde{i}}
       \leqslant \langle {\tilde{i}}_{\tilde Y},\tilde\lambda'_a\rangle\ \textrm{for any}\ {\tilde{i}}\in \tilde{I}, a=1,2.
\]
and so $x_1^+\xi_{-\tilde\lambda'_1}, x_2^-\eta_{\tilde\lambda'_2}\neq 0$. By \eqref{25.2.1}, $(x_1\diamondsuit_{\tilde\zeta'}x_2)c\neq 0$. Moreover, since the nonzero images of the canonical basis elements in $\dot{\tilde{\mathbf{B}}}$ under evaluation at $c$, are distinct canonical basis
elements in $\mathbf{B}({^{\omega}\Lambda_{\tilde{\lambda}'_1}\otimes \Lambda_{\tilde{\lambda}'_2}})$, a contradiction. Hence $\dot b \star \dot b'-(b_1\theta_N b'_1)\diamondsuit_{\tilde\zeta'}(b_2\theta_N b'_2)=0$.

It remains to prove \eqref{eq:star-leading-tensor} for any sufficiently large $N$. We need to calculate $$(\dot{b} \star \dot{b}')c=(b_1\diamondsuit_{\tilde{\zeta}} b_2)(\theta_N\diamondsuit_{\tilde{\zeta}} \theta_N)(b'_1\diamondsuit_{\tilde{\zeta}'} b'_2)c\in {^{\omega}\Lambda_{\tilde{\lambda}'_1}}\otimes \Lambda_{\tilde{\lambda}'_2}.$$
We write 
$$(b'_1\diamondsuit_{\zeta'}b'_2)(\xi_{-\lambda_1}\otimes \eta_{\lambda_2})=b_1'^+\xi_{-\lambda_1}\diamondsuit b_2'^-\eta_{\lambda_2}=\sum\pi_{x'y'}x'^+\xi_{-\lambda_1}\otimes y'^-\eta_{\lambda_2}$$
where the summation is taken over $x',y'\in\mathbf B$ such that $x'^+\xi_{-\lambda_1},y'^-\eta_{\lambda_2}\neq 0$. By \eqref{24.3.3}, we have
\begin{equation}\label{24.3.3'}
\pi_{b'_1,b'_2}=1;\ 
 \pi_{x',y'}\in v^{-1}\mathbb Z[v^{-1}]\ \textrm{if}\ (x',y')\neq (b'_1,b'_2)    
\end{equation}
Moreover, $\pi_{x',y'}\neq 0$ implies that $|x'|-|y'|=|\dot b'|$ (see \cite[Theorem 24.3.3 (b)]{Lusztig-1993}). Then by \eqref{characterizing condition}, 
\begin{equation}\label{transition matrix}
 (b'_1\diamondsuit_{\tilde\zeta'}b'_2)c=\iota((b'_1\diamondsuit_{\zeta'}b'_2)(\xi_{-\lambda_1}\otimes \eta_{\lambda_2}))
   =\sum \pi_{x',y'}\,x'^+\xi\otimes y'^-\eta,
\end{equation}
where the summation is taken over $x',y'\in\mathbf B$ such that $x'^+\xi,y'^-\eta\neq 0$.

For any $\mathbf n=(n_i)_{i\in I}$ with $0\leqslant n_i\leqslant N$ for any $i\in I$, we denote $\theta_{\mathbf n}=\prod_{i\in I}\theta_{i'}^{(n_i)}$. In particular, $\theta_{\mathbf N}=\theta_N$ for $\mathbf N=(N)_{i\in I}$. Applying Lemma \ref{lem:marked-powers} to \eqref{transition matrix}, we obtain 
$$(\theta_N\diamondsuit_{\tilde\zeta}\theta_N) (b'_1\diamondsuit_{\tilde\zeta'}b'_2)c
  =\sum_{x'^+\xi,y'^-\eta\neq 0} \pi_{x',y'}\sum_{\mathbf{0}\leqslant \mathbf t\leqslant \mathbf N}
      k_{\mathbf t}(x',y')\,
      (\theta_{\mathbf{N-t}}x')^+\xi
      \otimes (\theta_{\mathbf{N-t}}y')^-\eta,$$
where $k_{\mathbf0}(x',y')=1$ and $k_{\mathbf t}(x',y')\in v^{-\sum_{i\in I} t_i^2}\mathbb Z[v^{-1}]$. Then 
$$(\dot{b} \star \dot{b}')c=\sum_{x'^+\xi,y'^-\eta\neq 0} \pi_{x',y'}\sum_{\mathbf{0}\leqslant \mathbf t\leqslant \mathbf N} k_{\mathbf t}(x',y')\, (b_1\diamondsuit_{\tilde{\zeta}} b_2)\bigl((\theta_{\mathbf{N-t}}x')^+\xi \otimes (\theta_{\mathbf{N-t}}y')^-\eta\bigr).$$
We show that any term 
\begin{equation}\label{term}
k_{\mathbf t}(x',y')\, (b_1\diamondsuit_{\tilde{\zeta}} b_2)\bigl((\theta_{\mathbf{N-t}}x')^+\xi \otimes (\theta_{\mathbf{N-t}}y')^-\eta\bigr)\ \textrm{with}\  \mathbf t\neq\mathbf0  
\end{equation}
lies in $v^{-1}\mathcal L$ for sufficiently large $N$. By Lemma \ref{lem:marked-action}, there is
an integer $C\geqslant 0$ such that
\[
 (b_1\diamondsuit_{\tilde{\zeta}} b_2)\bigl((\theta_{\mathbf{N-t}}x')^+\xi \otimes (\theta_{\mathbf{N-t}}y')^-\eta\bigr)
 \in v^C\mathcal L
\]
for any $N$, $\mathbf 0<\mathbf t\leqslant \mathbf N$ and $x'^+\xi,y'^-\eta\neq 0 $ in summation. For those $\mathbf{t}$ such that $\sum_{i\in I} t_i^2>C$, we have $k_{\mathbf t}(x',y')\in v^{-C-1}\mathbb Z[v^{-1}]$.
Thus the corresponding term \eqref{term} lies in $v^{-1}\mathcal L$. It remains to consider finitely many $\mathbf{t}$ such that $0<\sum_{i\in I} t_i^2\leqslant C$. For any such $\mathbf t$, $N-t_i\rightarrow \infty$ when $N\rightarrow\infty$ for any $i\in I$. Applying the second assertion of Lemma~\ref{lem:marked-action} for $\mathbf{n}=\mathbf{N}-\mathbf{t}$, we obtain
\[
 (b_1\diamondsuit_{\tilde\zeta}b_2)\bigl((\theta_{\mathbf{N-t}}x')^+\xi \otimes (\theta_{\mathbf{N-t}}y')^-\eta\bigr)
 \in\mathcal L
\]
for any sufficiently large $N$, $\mathbf{t}$ such that $0<\sum_{i\in I} t_i^2\leqslant C$ and $x'^+\xi,y'^-\eta\neq 0$. Since $k_{\mathbf t}(x',y')\in v^{-1}\mathbb Z[v^{-1}]$ whenever
$\mathbf t\ne\mathbf0$, the corresponding terms $\eqref{term}$ also lie in $v^{-1}\mathcal L$.

For $\mathbf t=\mathbf0$, by Lemma \ref{lem:marked-powers}, $k_{\mathbf0}(x',y')=1$. Applying the second assertion of Lemma~\ref{lem:marked-action} for $\mathbf n=\mathbf N$, we obtain
\[
 \begin{aligned}
 (\dot b \star \dot b')c
 &\equiv \sum_{x'^+\xi,y'^-\eta\neq 0}\pi_{x',y'}\,
    (b_1\theta_N x')^+\xi\otimes(b_2\theta_N y')^-\eta \mod v^{-1}\mathcal L\\
 &\equiv (b_1\theta_N b'_1)^+\xi\otimes(b_2\theta_N b'_2)^-\eta \mod v^{-1}\mathcal L
 \end{aligned}
\]
for sufficiently large $N$. For the second congruence, each pure tensor $(b_1\theta_N x')^+\xi\otimes(b_2\theta_N y')^-\eta$ lies in
$\mathcal L$ by Lemma \ref{lem:sep}, while the coefficients satisfy
$\pi_{b'_1,b'_2}=1$ and $\pi_{x',y'}\in v^{-1}\mathbb Z[v^{-1}]$ for any other $x',y'$ by \eqref{24.3.3'}. This proves \eqref{eq:star-leading-tensor}.
\end{proof}

\subsection{Estimates for the action on tensor product}
\label{sec:action-estimates}

We prove the following Lemma \ref{lem:marked-powers} and Lemma \ref{lem:marked-action} to complete the proof of Proposition \ref{prop:product}. We retain the notation of the proof of Proposition~\ref{prop:product}. Fix $x',y'\in \mathbf{B}$ occurring in \eqref{transition matrix} so that $|x'|-|y'|=|\dot b'|$, otherwise the coefficient $\pi_{x',y'}$ is zero.

For any nonzero $f(v)=P(v)/Q(v)\in\mathbb Q(v)$, where
$P(v),Q(v)\in\mathbb Z[v]$, we define
\[
 \deg_\infty f(v)=\deg P(v)-\deg Q(v).
\]
By convention, we set $\deg_\infty 0=-\infty$. 

\begin{lemma}\label{lem:marked-powers}
We have 
\[
 (\theta_N\diamondsuit_{\tilde\zeta}\theta_N)
       (x'^+\xi\otimes y'^-\eta)
  =\sum_{\mathbf{0}\leqslant \mathbf t\leqslant \mathbf N}
      k_{\mathbf t}(x',y')\,
      (\theta_{\mathbf{N-t}}x')^+\xi
      \otimes (\theta_{\mathbf{N-t}}y')^-\eta,
\]
with coefficients 
$$k_{\mathbf t}(x',y')=\prod_{i\in I}v^{-t_i(t_i+r_i+|y'|_i)}\begin{bmatrix}
t_i+r_i+|y'|_i\\ t_i\end{bmatrix}.$$ 
In particular, $k_{\mathbf0}(x',y')=1$ and
$\deg_\infty k_{\mathbf t}(x',y')=-\sum_{i\in I}t_i^2$.
\end{lemma}
\begin{proof}
Since the new vertices $i'\in I'$ are pairwise orthogonal and
$\langle i'_{\tilde Y},\tilde\zeta\rangle=0$ for all $i\in I$, by
\cite[Proposition 25.3.2]{Lusztig-1993}, we have 
$\theta_N\diamondsuit_{\tilde\zeta}\theta_N
 =\theta_N^+\theta_N^-1_{\tilde\zeta}$.
Since $F_{i'}x'^+\xi=0$ for every $i\in I$, we have
\[
\theta_N^-(x'^+\xi\otimes y'^-\eta)
   =x'^+\xi\otimes(\theta_Ny')^-\eta.
\]

For any $i\in I$, by \cite[\S 3.1.5]{Lusztig-1993}, $\Delta(E_{i'}^{(N)})=\sum_{t_i=0}^N v^{t_i(N-t_i)} E_{i'}^{(N-t_i)}K_{i'_{\tilde{Y}}}^{t_i}\otimes E_{i'}^{(t_i)}$. For the first factor, by \eqref{lift lambda'_1,lambda'_2} and $|x'|-|y'|=|\dot{b}'|$, we have 
$$E_{i'}^{(N-t_i)}K_{i'_{\tilde{Y}}}^{t_i}x'^+\xi=v^{-t_i(N+r_i+|y'|_i)}E_{i'}^{(N-t_i)}x'^+\xi.$$ 
For the second factor, since $E_{i'}y'^-\eta=0$, by \cite[\S 3.4.2]{Lusztig-1993} and \eqref{lift lambda'_1,lambda'_2}, we have 
\begin{align*}
E_{i'}^{(t_i)}(\theta_Ny')^-\eta=\begin{bmatrix}
t_i+r_i+|y'|_i\\ t_i\end{bmatrix}(\prod_{j'\neq i'}F_{j'}^{(N)}F_{i'}^{(N-t_i)}y'^-)\eta.
\end{align*}
Repeating above process for any other $j'\neq i'$ in $I'$, we obtain the desired expression.

By definition, the quantum binomial coefficient
has degree $t_i(r_i+|y'|_i)$, and so the coefficient $k_{\mathbf t}(x',y')$ has degree $\sum_{i\in I}-t_i(t_i+r_i+|y'|_i)+t_i(r_i+|y'|_i)=-\sum_{i\in I}t_i^2$.
\end{proof}

It remains to calculate the action of $b_1\diamondsuit_{\tilde{\zeta}} b_2$ on each term $(\theta_{\mathbf{N-t}}x')^+\xi \otimes (\theta_{\mathbf{N-t}}y')^-\eta$ in the summation in Lemma \ref{lem:marked-powers}. 

\begin{lemma}\label{lem:marked-action}
There exists a positive integer $C$ such that
\[
 (b_1\diamondsuit_{\tilde{\zeta}} b_2)
 \bigl((\theta_{\mathbf{n}}x')^+\xi \otimes (\theta_{\mathbf{n}}y')^-\eta\bigr)
 \in v^C\mathcal L
\]
for every $N\in\mathbb N$, every $\mathbf n=(n_i)_{i\in I}$ with
$0\leqslant n_i\leqslant N$ for all $i\in I$, and every pair $(x',y')$
occurring in \eqref{transition matrix}. Moreover, if $\min\{n_i;i\in I\}$
is sufficiently large, then
\[
 (b_1\diamondsuit_{\tilde{\zeta}} b_2)
 \bigl((\theta_{\mathbf{n}}x')^+\xi \otimes (\theta_{\mathbf{n}}y')^-\eta\bigr)
 \equiv (b_1\theta_{\mathbf n}x')^+\xi
          \otimes(b_2\theta_{\mathbf n}y')^-\eta
 \pmod{v^{-1}\mathcal L},
\]
where the right-hand side is a pure tensor of canonical basis elements or zero.
\end{lemma}

\begin{proof}
We prove the lemma in three steps. We first bound the coefficients of in the linear combination of $\tilde{\mathbf{B}}$ for $x\theta_{\mathbf n}y\in \tilde{\mathbf{f}}$ with $x,y\in \mathbf{B}$, and show that it is a canonical basis element when all $n_i$ are sufficiently large. Then we use this bound to prove the first assertion, and show that if all $n_i$ are sufficiently large, then the terms in which a negative generator acts on $x'^+\xi$ or a positive generator acts on $y'^-\eta$ have strictly negative degree. Finally, to identify the sum of the remaining terms, we evaluate the same comultiplication on a tensor of a lowest weight vector and a highest weight vector with the same respective weights on the original vertices. The scalars in the calculation, and the terms of strictly negative degree just described vanish on this tensor.

\smallskip\noindent
\emph{Step 1. Multiplication in the negative part.}
We apply the norm calculation in the proof of Lemma \ref{lem:sep} with
$\theta_N$ replaced by $\theta_{\mathbf n}$. For fixed $x,y\in\mathbf B$, we have
\[
(x\theta_{\mathbf n}y,x\theta_{\mathbf n}y)=(\theta_{\mathbf n},\theta_{\mathbf n})
   \sum_{\nu\in \mathbb{N}[I]}v^{\nu\cdot\nu-2\sum_{i\in I} n_i\nu_i}A_\nu,
\]
where only finitely many $\nu$ occur, the $A_\nu$ are independent of
$\mathbf n$, and $A_0=(x,x)(y,y)$. Since the new vertices are pairwise
orthogonal, $\theta_{\mathbf n}\in\tilde{\mathbf B}$
by \cite[\S14.5.3 and Theorem~14.3.2(c)]{Lusztig-1993}. By almost
orthonormality, its norm has degree zero and leading coefficient~$1$.
Hence the degree of $(x\theta_{\mathbf n}y,x\theta_{\mathbf n}y)$ is
bounded above independently of $\mathbf n$. Write
$x\theta_{\mathbf n}y=\sum_{\tilde b\in\tilde{\mathbf B}}p_{\tilde b}\tilde b$.
If $x\theta_{\mathbf n}y\ne0$, then almost orthonormality gives
\[
 \deg_\infty(x\theta_{\mathbf n}y,x\theta_{\mathbf n}y)
   =2\max\{\deg_\infty p_{\tilde{b}}; \tilde{b}\in \tilde{\mathbf B}\}.
\]
The leading coefficient of the norm is the sum of the squares of the
leading coefficients of the $p_{\tilde b}$ of maximal degree. Thus
the degrees of the $p_{\tilde b}$ are bounded above independently of
$\mathbf n$. Moreover, if $\min_{i\in I}n_i$ is sufficiently large, each
term $v^{\nu\cdot\nu-2\sum_{i\in I} n_i\nu_i}A_\nu$ with $\nu\ne0$ has
strictly negative degree. By the same reason as in the proof of Lemma~\ref{lem:sep}, $x\theta_{\mathbf n}y\in\tilde{\mathbf B}$.

\smallskip\noindent
\emph{Step 2. Degree estimates.}
For any $\mathbf n=(n_i)$, let $-\lambda_{1,\mathbf n}$ and $\lambda_{2,\mathbf n}$ be the restrictions of the weights of $(\theta_{\mathbf{n}}x')^+\xi \in {^{\omega}\Lambda_{\tilde{\lambda}'_1}}$ and $(\theta_{\mathbf{n}}y')^-\eta\in \Lambda_{\tilde{\lambda}'_2}$ to the original coweight lattice $Y$. Here the restriction of a weight $\tilde{\lambda}\in \tilde{X}$ is the weight in $X$ corresponding to the function $Y\rightarrow \mathbb{Z},\mu\mapsto\langle \mu,\tilde{\lambda}\rangle$. Then we have for any $i\in I$,
\begin{align*}
&\langle i_Y,\lambda_{1,\mathbf n}\rangle=\langle i_Y,\lambda'_1-|x'|_X\rangle +n_i,\\
&\langle i_Y,\lambda_{2,\mathbf n}\rangle=\langle i_Y,\lambda'_2-|y'|_X\rangle +n_i 
\end{align*}
Since $|x'|-|y'|=|\dot{b}'|$, we have $\zeta=\lambda_{2,\mathbf n}-\lambda_{1,\mathbf n}$. By the triangular decomposition of $\dot{\mathbf{U}}$ and the definition of $\iota_{\zeta,\tilde{\zeta}}$, we can write 
$$b_1\diamondsuit_{\zeta} b_2=\sum_{\mathbf i,\mathbf{j}}c_{\mathbf i,\mathbf{j}}E_{\mathbf i}F_{\mathbf j}1_{\zeta},\ b_1\diamondsuit_{\tilde{\zeta}} b_2=\sum_{\mathbf i,\mathbf{j}}c_{\mathbf i,\mathbf{j}}E_{\mathbf i}F_{\mathbf j}1_{\tilde{\zeta}}$$
where $\mathbf{i},\mathbf{j}$ are sequences in $I$ and $E_{\mathbf i},F_{\mathbf j}$ are corresponding monomials in Chevalley generators. The coefficients $c_{\mathbf i,\mathbf{j}}$ are independent of $\mathbf n$.

Consider the action of $E_{\mathbf i}F_{\mathbf j}$ on $(\theta_{\mathbf{n}}x')^+\xi\otimes (\theta_{\mathbf{n}}y')^-\eta$. By \cite[\S 3.1.5 \& Lemma 1.2.11]{Lusztig-1993}, we have 
\begin{align*}
&\Delta(E_{\mathbf i})=\sum v^{\kappa(\mathbf i',\mathbf i'')}E_{\mathbf i'}K_{|E_{\mathbf i''}|_Y}\otimes E_{\mathbf i''},\\
&\Delta(F_{\mathbf j})=\sum v^{\kappa'(\mathbf j',\mathbf j'')} F_{\mathbf j'}\otimes F_{\mathbf j''}K_{-|F_{\mathbf j'}|_Y},    
\end{align*} 
The sums run over subsequences $\mathbf i',\mathbf j'$ and their complementary subsequences $\mathbf i'',\mathbf j''$, with positions counted separately. The powers $\kappa(\mathbf i',\mathbf i''),\kappa'(\mathbf j',\mathbf j'')$ are independent of $\mathbf n$. The degrees $|E_{\mathbf i''}|,|F_{\mathbf j'}|\in \mathbb{N}[I]$ are determined by $E_{\mathbf i''}\in \mathbf{U}(|E_{\mathbf i''}|),|F_{\mathbf j'}|\in \mathbf{U}(-|F_{\mathbf j'}|)$. Then $E_{\mathbf i}F_{\mathbf j}((\theta_{\mathbf{n}}x')^+\xi\otimes (\theta_{\mathbf{n}}y')^-\eta)$ is equal to 
\begin{align*}
\sum v^{C_{\mathbf i',\mathbf i'',\mathbf j',\mathbf j''}-\langle |E_{\mathbf i''}|_Y,\lambda_{1,\mathbf{n}}\rangle-\langle |F_{\mathbf{j'}}|_Y,\lambda_{2,\mathbf{n}}\rangle}E_{\mathbf i'}F_{\mathbf j'} (\theta_{\mathbf{n}}x')^+\xi\otimes E_{\mathbf i''}F_{\mathbf j''}(\theta_{\mathbf{n}}y')^-\eta,
\end{align*}
where $C_{\mathbf i',\mathbf i'',\mathbf j',\mathbf j''}$ is independent of $\mathbf{n}$. We use the commutation relation $E_iF_j-F_jE_i=\delta_{i,j}\frac{K_{i_Y}-K_{-i_Y}}{v-v^{-1}}$ to write the second factor 
$$E_{\mathbf i''}F_{\mathbf j''}(\theta_{\mathbf{n}}y')^-\eta=\sum d_{\mathbf{i'''},\mathbf{j'''}} F_{\mathbf j'''}E_{\mathbf i'''}(\theta_{\mathbf{n}}y')^-\eta,$$
where $\mathbf{i}''',\mathbf{j}'''$ are subsequences of $\mathbf i'',\mathbf{j}''$ respectively. Before collecting like terms, let $k_1,\ldots,k_r\in I$ be the vertices at which the relation is used in a given summand. Each such use deletes one $E_{k_s}$ and one $F_{k_s}$. The contribution to the coefficient is $\prod_{s=1}^r [\langle {k_{s}}_Y,\lambda_{2,\mathbf{n}}+\rho_s\rangle]$, where $\rho_s\in X$ is independent of $\mathbf{n}$. Hence $E_{\mathbf i}F_{\mathbf j}((\theta_{\mathbf{n}}x')^+\xi\otimes (\theta_{\mathbf{n}}y')^-\eta)$ is equal to 
$$\sum c_{\mathbf i',\mathbf i''',\mathbf j',\mathbf j'''} E_{\mathbf i'}F_{\mathbf j'} (\theta_{\mathbf{n}}x')^+\xi\otimes F_{\mathbf j'''}E_{\mathbf i'''}(\theta_{\mathbf{n}}y')^-\eta.$$
For $m\ne0$, we have $\deg_\infty[m]=|m|-1$, and $\deg_\infty 0=-\infty$. Since $\langle i_Y,\lambda_{1,\mathbf n}\rangle-n_i,\langle i_Y,\lambda_{2,\mathbf n}\rangle-n_i$ are independent of $\mathbf n$ for any $i\in I$, the contribution of each summand to $c_{\mathbf i',\mathbf i''',\mathbf j',\mathbf j'''}$, denoted by $c$, satisfies
\begin{align*}
\deg_{\infty} c&\leqslant C_{\mathbf i',\mathbf i'',\mathbf j',\mathbf j''}-\langle |E_{\mathbf i''}|_Y,\lambda_{1,\mathbf{n}}\rangle-\langle |F_{\mathbf{j'}}|_Y,\lambda_{2,\mathbf{n}}\rangle+\sum_{s=1}^r(\lvert\langle {k_{s}}_Y,\lambda_{2,\mathbf{n}}+\rho_s\rangle\rvert-1)\\
&\leqslant C'_{\mathbf i',\mathbf i''',\mathbf j',\mathbf j'''}-\sum_{i\in I}n_i(|E_{\mathbf i'''}|_i+|F_{\mathbf{j}'}|_i),
\end{align*}
where $C'_{\mathbf i',\mathbf i''',\mathbf j',\mathbf j'''}$ is independent of $N$ and $\mathbf{n}$. The last inequality uses $|n_i+a|\leqslant n_i+|a|$ and $|E_{\mathbf i''}|-|E_{\mathbf i'''}|=\sum_{s=1}^r k_s$. Hence  $(b_1\diamondsuit_{\tilde{\zeta}}b_2)((\theta_{\mathbf{n}}x')^+\xi\otimes (\theta_{\mathbf{n}}y')^-\eta)$ can be written as
\begin{equation}\label{action result}
\sum_{\mathbf p,\mathbf q,\mathbf r,\mathbf s}a_{\mathbf p,\mathbf q,\mathbf r,\mathbf s}(\mathbf n)E_{\mathbf p}F_{\mathbf q} (\theta_{\mathbf{n}}x')^+\xi\otimes F_{\mathbf r}E_{\mathbf s}(\theta_{\mathbf{n}}y')^-\eta, 
\end{equation}
where $\mathbf p,\mathbf q,\mathbf r,\mathbf s$ are sequences in $I$, and the degrees of the coefficients satisfy 
\begin{equation}\label{bound}
\deg_{\infty} a_{\mathbf p,\mathbf q,\mathbf r,\mathbf s}(\mathbf n)\leqslant C-\sum_{i\in I}n_i(|\mathbf q|_i+|\mathbf s|_i),    
\end{equation}
where $C$ is independent of $N$ and $\mathbf n$. A single bound suffices since only finitely many $\mathbf{i},\mathbf{j}$ and straightening choices $[E_i,F_j]$ occur. For each term 
$$E_{\mathbf p}F_{\mathbf q}(\theta_{\mathbf{n}}x')^+\xi\otimes F_{\mathbf r}E_{\mathbf s}(\theta_{\mathbf{n}}y')^-\eta=E_{\mathbf p}\theta_{\mathbf{n}}^+F_{\mathbf q}x'^+\xi\otimes F_{\mathbf r}\theta_{\mathbf{n}}^-E_{\mathbf s}y'^-\eta$$ 
in the summation \eqref{action result}, using the inclusions $\iota_1,\iota_2$ in the proof of Lemma \ref{lift}, we can express $F_{\mathbf q}x'^+\xi$ into the linear combination of $x_2^+\xi$ with $x_2\in \mathbf{B}$, and express $E_{\mathbf s}y'^-\eta$ in the linear combination of $y_2^-\eta$ with $y_2\in \mathbf{B}$. Then expressing the monomials $E_{\mathbf p},F_{\mathbf r}$ into the linear combinations of $x_1^+,y_1^-$ with $x_1,y_1\in \mathbf{B}$ respectively gives
$$E_{\mathbf p}F_{\mathbf q}(\theta_{\mathbf{n}}x')^+\xi\otimes F_{\mathbf r}E_{\mathbf s}(\theta_{\mathbf{n}}y')^-\eta=\sum c_{x_1,x_2,y_1,y_2} (x_1\theta_{\mathbf{n}}x_2)^+\xi\otimes (y_1\theta_{\mathbf{n}}y_2)^-\eta,$$
where $x_1,x_2,y_1,y_2\in \mathbf{B}$ and the coefficients $c_{x_1,x_2,y_1,y_2}$ are independent of $\mathbf n$. By the norm calculation in Step 1, when we write the products $x_1\theta_{\mathbf{n}}x_2,y_1\theta_{\mathbf{n}}y_2$ into the linear combinations of $\tilde{\mathbf{B}}$, the degrees of the coefficients are uniformly bounded independently of $\mathbf{n}$. Together with \eqref{bound} and the integrality of the action, this proves the first assertion. This also shows that if $\min\{n_i;i\in I\}$ is sufficiently large, then each term with $(\mathbf q,\mathbf s)\ne(\varnothing,\varnothing)$ has coefficients of strictly negative degree with respect to the basis $\mathbf{B}({^{\omega}\Lambda_{\tilde\lambda'_1}})\otimes \mathbf{B}(\Lambda_{\tilde\lambda'_2})$.

\smallskip\noindent
\emph{Step 3. Comparison with leading terms.}
It remains to identify the sum of the terms with $\mathbf q=\mathbf s=\varnothing$. If $\min \{n_i;i\in I\}$ is sufficiently large, then $\lambda_{1,\mathbf n},\lambda_{2,\mathbf n}$ are dominant. Consider the action of $b_1\diamondsuit_{\zeta}b_2$ on $\xi_{-\lambda_{1,\mathbf n}}\otimes \eta_{\lambda_{2,\mathbf n}}$. By the same calculation as \eqref{action result}, $(b_1\diamondsuit_{\zeta}b_2)(\xi_{-\lambda_{1,\mathbf n}}\otimes \eta_{\lambda_{2,\mathbf n}})$ is equal to 
\begin{align*}
\sum_{\mathbf p,\mathbf q,\mathbf r,\mathbf s}a_{\mathbf p,\mathbf q,\mathbf r,\mathbf s}(\mathbf n)E_{\mathbf p}F_{\mathbf q} \xi_{-\lambda_{1,\mathbf n}}\otimes F_{\mathbf r}E_{\mathbf s} \eta_{\lambda_{2,\mathbf n}}
=\sum_{\mathbf p,\mathbf r}a_{\mathbf p,\varnothing,\mathbf r,\varnothing}(\mathbf n) E_{\mathbf p}\xi_{-\lambda_{1,\mathbf n}} \otimes F_{\mathbf r}\eta_{\lambda_{2,\mathbf n}}
\end{align*}
with the same coefficients as in \eqref{action result}, since $-\lambda_{1,\mathbf n}$ and $\lambda_{2,\mathbf n}$ are the restrictions of the weights of $(\theta_{\mathbf{n}}x')^+\xi \in {^{\omega}\Lambda_{\tilde{\lambda}'_1}}$ and $(\theta_{\mathbf{n}}y')^-\eta\in \Lambda_{\tilde{\lambda}'_2}$ respectively. The other terms vanish because $F_{\mathbf q}\xi_{-\lambda_{1,\mathbf n}}=0$ for $\mathbf q\ne\varnothing$ and $E_{\mathbf s}\eta_{\lambda_{2,\mathbf n}}=0$ for $\mathbf s\ne\varnothing$. Let
\[
 G_{\mathbf n}
   =\sum_{\mathbf p,\mathbf r}a_{\mathbf p,\varnothing,\mathbf r,\varnothing}(\mathbf n) \theta_{\mathbf p} \otimes \theta_{\mathbf r}\in\mathbf f\otimes\mathbf f .
\]
Then the map 
$$\mathbf{f}\otimes\mathbf f\rightarrow {^{\omega}\Lambda_{\lambda_{1,\mathbf n}}}\otimes \Lambda_{\lambda_{2,\mathbf n}},\ x\otimes y\mapsto x^+\xi_{-\lambda_{1,\mathbf n}}\otimes y^-\eta_{\lambda_{2,\mathbf n}}$$
sends $G_{\mathbf n}$ to $(b_1\diamondsuit_{\zeta}b_2)(\xi_{-\lambda_{1,\mathbf n}}\otimes \eta_{\lambda_{2,\mathbf n}})$. By \eqref{25.2.1} and \eqref{24.3.3}, we have
\begin{equation}\label{G}
(b_1\diamondsuit_{\zeta}b_2)(\xi_{-\lambda_{1,\mathbf n}}\otimes \eta_{\lambda_{2,\mathbf n}})\equiv b_1^+\xi_{-\lambda_{1,\mathbf n}}\otimes b_2^-\eta_{\lambda_{2,\mathbf n}}\mod v^{-1}\mathbb{Z}[v^{-1}][\mathbf{B}({^{\omega}\Lambda_{\lambda_{1,\mathbf n}}})\otimes \mathbf{B}(\Lambda_{\lambda_{2,\mathbf n}})].   
\end{equation}
The degrees $|\theta_{\mathbf p}|,|\theta_{\mathbf r}|$ occurring in $G_{\mathbf n}$ are bounded independently of $\mathbf{n}$. Increasing $\min\{n_i;i\in I\}$ if necessary, the maps $\mathbf{f}\rightarrow {^{\omega}\Lambda_{\lambda_{1,\mathbf n}}},x\mapsto x^+\xi_{-\lambda_{1,\mathbf n}}$ and $\mathbf{f}\rightarrow \Lambda_{\lambda_{2,\mathbf n}},y\mapsto y^-\eta_{\lambda_{2,\mathbf n}}$ send distinct canonical basis elements of $\mathbf{f}$ of degrees occurring in $G_{\mathbf n}$ or in the respective degrees $|b_1|,|b_2|$ to distinct nonzero canonical basis elements of ${^{\omega}\Lambda_{\lambda_{1,\mathbf n}}}$ and $\Lambda_{\lambda_{2,\mathbf n}}$ respectively. So \eqref{G} gives 
\[
 G_{\mathbf n}\equiv b_1\otimes b_2\mod v^{-1}\mathbb Z[v^{-1}][\mathbf B\otimes\mathbf B].
\]
Apply to this congruence the map
\[
 \mathbf{f}\otimes\mathbf f\rightarrow {^{\omega}\Lambda_{\tilde{\lambda}'_1}}\otimes \Lambda_{\tilde{\lambda}'_2}, x\otimes y\longmapsto
 (x\theta_{\mathbf n}x')^+\xi\otimes
 (y\theta_{\mathbf n}y')^-\eta.
\]
By Step 1, when $\min\{n_i;i\in I\}$ is sufficiently large, this map sends each pure tensor of canonical basis elements in the
fixed degrees under consideration to a pure tensor of canonical basis elements or to zero. Its value on $G_{\mathbf n}$ is 
$$\sum_{\mathbf p,\mathbf r}a_{\mathbf p,\varnothing,\mathbf r,\varnothing}(\mathbf n)E_{\mathbf p} (\theta_{\mathbf{n}}x')^+\xi\otimes F_{\mathbf r}(\theta_{\mathbf{n}}y')^-\eta.$$
Comparing with \eqref{action result}, the remaining terms satisfy $(\mathbf q,\mathbf s)\ne(\varnothing,\varnothing)$ and have coefficients of strictly negative degrees which are uniformly bounded by Step 2.

Finally, $b_1\diamondsuit_{\tilde{\zeta}}b_2$ preserves the integral form of the tensor
product, and both $(\theta_{\mathbf{n}}x')^+\xi \otimes (\theta_{\mathbf{n}}y')^-\eta$ and
the tensor on the right of the asserted congruence are integral.
Thus every coefficient of the difference belongs to
$\mathbb Z[v,v^{-1}]$ and has strictly negative degree. It
therefore belongs to $v^{-1}\mathbb Z[v^{-1}]$, as required.
\end{proof}

\subsection{Positivity of multiplication and tensor product actions}

We now deduce Lusztig's positivity conjecture for multiplication
from the corresponding result for comultiplication.

\begin{corollary}\label{cor:positive}
Let $\dot\U$ be the modified quantum group associated with any root
datum of symmetric Cartan datum. For all $\dot b,\dot b'\in\dot\B$,
\[
 \dot b\dot b'\in\mathbb N[v,v^{-1}][\dot\B].
\]
\end{corollary}

\begin{proof}
Consider first the simply connected root datum. We may assume that
the weights of $\dot b$ and $\dot b'$ match, since otherwise $\dot b \dot b'=0$. Choose $N\in \mathbb{N}$ sufficiently large that
$\dot b\star\dot b'$ belongs to $\dot{\tilde{\mathbf B}}$, as in
Proposition~\ref{prop:product}. Theorem~\ref{thm:main} then identifies
each coefficient of $\dot b\dot b'$ with a coefficient in the
comultiplication of this canonical basis element. Since the enlarged Cartan
datum is symmetric, every such coefficient belongs to
$\mathbb N[v,v^{-1}]$ by \cite[Theorem~6.3]{Fang-He-2026}.

For an arbitrary root datum of the same Cartan type, the comparison
in \cite[\S6.6]{Fang-He-2026} identifies the multiplication
coefficients for matching weights with those for the simply
connected root datum. This proves the assertion.
\end{proof}

We next consider the action on tensor products. Let $m,n\geqslant 0$
with $m+n\geqslant 1$, and let
$\boldsymbol\lambda=(\lambda_1,\ldots,\lambda_{m+n})$ be a sequence
of dominant weights. Write
\[
 \Lambda(\boldsymbol\lambda)
 =
 {}^\omega\Lambda_{\lambda_1}\otimes\cdots\otimes
 {}^\omega\Lambda_{\lambda_m}\otimes
 \Lambda_{\lambda_{m+1}}\otimes\cdots\otimes\Lambda_{\lambda_{m+n}},
\]
with canonical basis $\B(\Lambda(\boldsymbol\lambda))$ as in
\cite[\S7.1]{Fang-He-2026}. When $m=0$ or $n=0$, positivity of the
action was proved in \cite[Theorem~7.4(2)]{Fang-He-2026}. When both
$m$ and $n$ are positive, the following result removes the spherical
parabolic hypothesis on $\dot b$.

\begin{corollary}\label{cor:tensor-action}
For any root datum of symmetric Cartan type and any $\dot b\in\dot\B$,
\[
 \dot b\bigl(\B(\Lambda(\boldsymbol\lambda))\bigr)
 \subseteq
 \mathbb N[v,v^{-1}]\bigl[\B(\Lambda(\boldsymbol\lambda))\bigr].
\]
\end{corollary}

\begin{proof}
Consider first ${}^\omega\Lambda_\lambda\otimes\Lambda_\mu$, and put
$c=\xi_{-\lambda}\otimes\eta_\mu$. By
\cite[Theorem~1.4]{Fang-He-2026}, its canonical basis consists of the
nonzero vectors $\dot a c$ with
$\dot a\in\dot\B\cap\dot\U 1_{\mu-\lambda}$. Since
\[
 \dot b(\dot a c)=(\dot b\dot a)c,
\]
Corollary~\ref{cor:positive} and \eqref{25.2.1} gives the assertion for this tensor product.

For an arbitrary number of factors, the proof of
\cite[Theorem~7.4(2)]{Fang-He-2026} realizes
$\Lambda(\boldsymbol\lambda)$ as a based submodule of a tensor
product of one lowest weight module and one highest weight module
for an enlarged symmetric Cartan datum. Applying the two-factor
case there and restricting to this based submodule gives the result.
\end{proof}

\section{Examples}\label{sec:examples}
\subsection{Examples in type $A_1$}\label{subsec:rank-one-match}

Let $I=\{i\}$ with $i\cdot i=2$. For the simply connected root datum
we have $Y=X=\mathbb Z$, $i_Y=1$ and $i_X=2$. The canonical basis of
$\mathbf f$ consists of the divided powers $\theta_i^{(n)}$,
$n\in\mathbb N$. By \cite[Proposition~25.3.2]{Lusztig-1993}, the
canonical basis elements of $\dot{\mathbf U}$ consists of
\[
 \theta_i^{(k)}\diamondsuit_t\theta_i^{(l)}
   =\begin{cases}
       E_i^{(k)}F_i^{(l)}1_t,& t\leqslant l-k,\\
       F_i^{(l)}E_i^{(k)}1_t,& t\geqslant l-k,
     \end{cases}
 \qquad (t\in\mathbb Z,\ k,l\in\mathbb N),
\]
where the two expressions agree for $t=l-k$.

The thickened diagram is $i\! -\! i'$ of type $A_2$. The canonical basis of $\tilde{\mathbf f}$ consists of the elements
\[
 \theta_i^{(p)}\theta_{i'}^{(N)}\theta_i^{(r)},\qquad
 \theta_{i'}^{(r)}\theta_i^{(N)}\theta_{i'}^{(p)}
 \qquad (N\geqslant p+r),
\]
with the two expressions identified when $N=p+r$
\cite[\S~14.5.4]{Lusztig-1993}. Thus $N\geqslant p+r$ suffices in
Lemma~\ref{lem:sep} for the pair $(\theta_i^{(p)},\theta_i^{(r)})$.
In both examples below, let
\[
 H_N=\theta_{i'}^{(N)}\diamondsuit_0\theta_{i'}^{(N)}
    =E_{i'}^{(N)}F_{i'}^{(N)}1_0.
\]

\smallskip\noindent
(1) Let $\dot b=E_i^{(2)}1_0$ and $\dot b'=F_i^{(2)}1_4$. Then \cite[\S~23.1.3]{Lusztig-1993} gives
\begin{equation}\label{eq:rank-one-divided}
 \dot b\dot b'=F_i^{(2)}E_i^{(2)}1_4
        +[4]F_iE_i1_4+\genfrac{[}{]}{0pt}{}{4}{2}1_4 .
\end{equation}
Take $N=4$. By \cite[Theorem~3.1]{Cui-2018}, we have
\[
 \dot b\star\dot b'
   =E_i^{(2)}E_{i'}^{(4)}F_{i'}^{(4)}F_i^{(2)}1_{\tilde{\zeta}'}
   =(\theta_i^{(2)}\theta_{i'}^{(4)})\diamondsuit_{\tilde{\zeta}'}
     (\theta_{i'}^{(4)}\theta_i^{(2)}),
\]
where $\tilde{\zeta}'=(4,-2)$. 
The terms in $\dot\Delta_{\tilde\zeta',0}(\dot b\star\dot b')$
with second factor $H_4$ is
\[
 \left(F_i^{(2)}E_i^{(2)}1_{\tilde\zeta'}
       +[4]F_iE_i1_{\tilde\zeta'}
       +\genfrac{[}{]}{0pt}{}{4}{2}1_{\tilde\zeta'}\right)
 \otimes H_4 .
\]
Thus all three coefficients in \eqref{eq:rank-one-divided} occur.

\smallskip\noindent
(2) Let $\dot b=E_iF_i1_{-2}$. Then $F_iE_i1_{-4}=(E_iF_i+[4])1_{-4}$ gives
\begin{equation}\label{eq:mixed-square}
 \dot b^2=[2]^2 E_i^{(2)}F_i^{(2)}1_{-2}+[4]E_iF_i1_{-2} .
\end{equation}
Take $N=4$. We denote $x=\theta_i\theta_{i'}^{(4)}\theta_i$, then 
\[
 \dot b\star\dot b
   =E_iF_iE_{i'}^{(4)}F_{i'}^{(4)}E_iF_i1_{\tilde{\zeta}'}
   =x^+x^-1_{\tilde{\zeta}'}=x\diamondsuit_{\tilde{\zeta}'}x,
\]
where $\tilde{\zeta}=(-2,0)$. For the second equality, we commute the generators at distinct vertices $F_iE_{i'}^{(4)}F_{i'}^{(4)}E_i=E_{i'}^{(4)}F_iE_iF_{i'}^{(4)}$, then $F_iE_i$ acts on $F_{i'}^{(4)}F_i1_{\tilde{\zeta}'}$ at $i$-weight $-2+4-2=0$, so it equals
$E_iF_i$ there. The last equality follows from
\cite[Theorem~3.1]{Cui-2018}. The terms in the coproduct with the second factor $H_4$ is
\[
 \dot{\Delta}_{\tilde{\zeta}',0}
           (\dot b\star\dot b)
     =([2]^2 E_i^{(2)}F_i^{(2)}1_{\tilde{\zeta}'}+[4]E_iF_i1_{\tilde{\zeta}'})\otimes H_4,
\]
Thus all coefficients in in \eqref{eq:mixed-square} occur.

\subsection{Example in type $A_2$}\label{subsec:a2-a4}

Let $i \! -\! j$ be the diagram of type $A_2$.
In the order $(i,j)$, let $\zeta=(-2,2)$ and $\zeta'=(4,-4)$.
Consider the canonical basis elements
\[
 \dot b=E_i^{(2)}F_j^{(2)}1_\zeta
       =\theta_i^{(2)}\diamondsuit_\zeta\theta_j^{(2)},\qquad
 \dot b'=E_j^{(2)}F_i^{(2)}1_{\zeta'}
       =\theta_j^{(2)}\diamondsuit_{\zeta'}\theta_i^{(2)}
\]
of $\dot{\mathbf U}$ (see \cite[Theorem~3.1]{Cui-2018}).
Their product is equal to
\begin{equation}\label{eq:a2-product}
\begin{aligned}
 \dot b\dot b'={}&
 (\theta_i^{(2)}\theta_j^{(2)})\diamondsuit_{\zeta'}
          (\theta_j^{(2)}\theta_i^{(2)})
 +\theta_i^{(2)}\diamondsuit_{\zeta'}\theta_i^{(2)}
 +\theta_j^{(2)}\diamondsuit_{\zeta'}\theta_j^{(2)}\\
 &+[2]\bigl((\theta_i^{(2)}\theta_j)\diamondsuit_{\zeta'}
                 (\theta_j\theta_i^{(2)})
           +(\theta_i\theta_j^{(2)})\diamondsuit_{\zeta'}
                 (\theta_j^{(2)}\theta_i)\bigr)\\
 &+[3]!\bigl((\theta_i\theta_j)\diamondsuit_{\zeta'}(\theta_j\theta_i)
           +\theta_i\diamondsuit_{\zeta'}\theta_i
           +\theta_j\diamondsuit_{\zeta'}\theta_j\bigr)
   +\genfrac{[}{]}{0pt}{}{4}{2}1_{\zeta'} .
\end{aligned}
\end{equation}

The thickened diagram is $i'\! -\! i\! -\! j\! -\! j'$ of type $A_4$.
We use this order for the coordinates of weights. Take $N=2$. We denote  $x=\theta_i^{(2)}\theta_{i'}^{(2)}\theta_{j'}^{(2)}\theta_j^{(2)}$ and
$y=\theta_j^{(2)}\theta_{i'}^{(2)}\theta_{j'}^{(2)}\theta_i^{(2)}$, then
\[
\dot b\star\dot b'=E_i^{(2)}F_j^{(2)}E_{i'}^{(2)}E_{j'}^{(2)}F_{i'}^{(2)}F_{j'}^{(2)}E_j^{(2)}F_i^{(2)}1_{\tilde\zeta'}
   =x^+y^-1_{\tilde\zeta'},
\]
where $\tilde\zeta'=(-2,4,-4,2)$. For the last equality, we commute the generators at distinct vertices $F_j^{(2)}E_{i'}^{(2)}E_{j'}^{(2)}F_{i'}^{(2)}F_{j'}^{(2)}E_j^{(2)}=E_{i'}^{(2)}E_{j'}^{(2)}F_j^{(2)}E_j^{(2)}F_{i'}^{(2)}F_{j'}^{(2)}$, then $F_j^{(2)}E_j^{(2)}$ acts on $F_{i'}^{(2)}F_{j'}^{(2)}F_i^{(2)}1_{\tilde{\zeta}'}$ at $j$-weight $-4+2+2=0$, so it equals to $E_j^{(2)}F_j^{(2)}$ there. In this example, $\dot b\star\dot b'=x\diamondsuit_{\tilde\zeta'}y$ is a canonical basis element already for $N=2$. Let $e_k=\Delta(E_k^{(2)})$ and $f_k=\Delta(F_k^{(2)})$. Then a direct computation in type $A_4$ gives
\begin{equation}\label{eq:a4-coproduct}
\begin{aligned}
\dot{\Delta}_{\tilde\zeta',0}(\dot b\star\dot b')
 &=e_i e_{i'}e_{j'}e_j\, f_j f_{i'}f_{j'}f_i
                    (1_{\tilde\zeta'}\otimes1_0)\\
 &=Z\otimes H_2
   +(\dot{b}\star \dot{b}')\otimes1_0
   +\sum Z_{\dot{b}''}\otimes \dot{b}'' ,
\end{aligned}
\end{equation}
where $\dot{b}''\in (\dot{\tilde{\mathbf{B}}}\cap\dot{\tilde{\mathbf{U}}}1_0)\setminus \{H_2,1_0\}$, and 
\begin{align*}Z&=\theta_i^{(2)}\theta_j^{(2)}\diamondsuit_{\tilde{\zeta}}\theta_j^{(2)}\theta_i^{(2)}+\theta_i^{(2)}\diamondsuit_{\tilde{\zeta}}\theta_i^{(2)}+\theta_j^{(2)}\diamondsuit_{\tilde{\zeta}}\theta_j^{(2)}\\
&+[2](\theta_i^{(2)}\theta_j\diamondsuit_{\tilde{\zeta}}\theta_j\theta_i^{(2)}+\theta_i\theta_j^{(2)}\diamondsuit_{\tilde{\zeta}}\theta_j^{(2)}\theta_i)\\
&+[3]!(\theta_i\theta_j\diamondsuit_{\tilde{\zeta}}\theta_j\theta_i+\theta_i\diamondsuit_{\tilde{\zeta}}\theta_i+\theta_j\diamondsuit_{\tilde{\zeta}}\theta_j)+\genfrac{[}{]}{0pt}{}{4}{2}1_{\tilde{\zeta}}.
\end{align*}

There are $29211$ nonzero terms in the canonical basis expansion in \eqref{eq:a4-coproduct}, exactly nine terms have the second factor $H_2$, and all coefficients in \eqref{eq:a2-product} occur term for term. For instance, the
coefficients of $(\theta_i\diamondsuit_{\tilde\zeta'}\theta_i)\otimes H_2$ and $1_{\tilde\zeta'}\otimes H_2$ are
\[
 [3]!=v^3+2v+2v^{-1}+v^{-3},\qquad 
 \genfrac{[}{]}{0pt}{}{4}{2}=v^4+v^2+2+v^{-2}+v^{-4}
\]
respectively. Both of them specialize to $6$ at $v=1$, and we can only distinguish them from each other in the quantization case.

The two expansions can be computed independently in the
\href{https://quantum-positivity-lab.xuhuahe.chatgpt.site/\#a2-a4-divided}
{\emph{Quantum Positivity Lab}}. To reproduce \eqref{eq:a4-coproduct}
with the linked example, set the second right weight to $(0,0,0,0)$.

\subsection{An example where the converse fails}\label{subsec:converse}

Theorem~\ref{thm:main} realizes each structure constant of
multiplication in $\dot{\mathbf U}$ with respect to $\dot{\mathbf B}$
as a structure constant of comultiplication in $\dot{\tilde{\mathbf U}}$
with respect to $\dot{\tilde{\mathbf B}}$.
Conversely, does every canonical basis coefficient of comultiplication occur as a
multiplication coefficient in some modified quantum group?
The answer is negative already in type $A_1$.

Let $\dot{\mathbf U}$ be a modified quantum group with canonical basis
$\dot{\mathbf B}$. For $\dot b,\dot b'\in\dot{\mathbf B}$, write
\[
 \dot b\dot b'=\sum_{\dot b''\in\dot{\mathbf B}}
      \dot m_{\dot b,\dot b'}^{\dot b''}\dot b''.
\]
Since the bar involution on $\dot{\mathbf U}$ preserves multiplication
and fixes each canonical basis element
(see \cite[Theorem~25.2.1]{Lusztig-1993}), we have
$\overline{\dot m_{\dot b,\dot b'}^{\dot b''}}
 =\dot m_{\dot b,\dot b'}^{\dot b''}$.
In $\dot{\mathbf U}(\mathfrak{sl}_2)$, we have
\[
 \dot\Delta_{0,0}(E_i^{(2)}1_0)
   = E_i^{(2)}1_0\otimes 1_0
     +v E_i1_0\otimes E_i1_0
     +1_0\otimes E_i^{(2)}1_0 .
\]
The coefficient of $E_i1_0\otimes E_i1_0$ is $v$.
Since $\bar v=v^{-1}\ne v$, it cannot be a multiplication coefficient.

\smallskip

A separate question is whether every multiplication structure constant
for the canonical basis of $\dot{\mathbf U}$ occurs as a multiplication
structure constant for the canonical basis of the negative part of a
larger quantum group. For symmetric Cartan data, this holds whenever
one factor is spherical parabolic; see
\cite[Theorem~4.3 and the proof of Theorem~6.6]{Fang-He-2026}.

\section*{AI declaration}
We used ChatGPT to develop the
\href{https://quantum-positivity-lab.xuhuahe.chatgpt.site/}%
{\textcolor{blue!60!black}{\underline{Quantum Positivity Lab}}}, an open-access tool for computing structure constants of multiplication and comultiplication
in modified quantum groups with respect to Lusztig's canonical bases.
These computations, including the comparison of multiplication in type
$A_2$ with comultiplication in type $A_4$ discussed in
\S\ref{subsec:a2-a4}, helped us discover the identity between the two
families of structure constants proved in this paper. We also used
ChatGPT for English-language editing, including improvements to grammar,
wording, and readability. We reviewed all suggested changes and take
full responsibility for the content of the paper.

\end{document}